\documentclass{tac}
\usepackage{xy}
\usepackage{MnSymbol}
\xyoption{all}

\def\mathcaldef#1{\expandafter\def\csname#1\endcsname{{\cal#1}}}

\def\v{``}

\def\qq{\quad\quad}
\def\qv{\qq ;\qq}
\def\iso{\,\cong\,}

\def\imp{\Rightarrow}
\def\la{\langle}
\def\ra{\rangle}
\def\adj{\dashv}
\def\op{^{\rm op}}
\def\ast{^*}
\def\st{^{-1}}
\def\ex{\exists}
\def\fa{\forall}
\def\comp{\lefthalfcap}
\def\tr{\triangleright}
\def\ov{\overline}

\def\tm{\times}
\def\otm{\otimes}
\def\dm{\Diamond}
\def\sq{\square}
\def\cd{\cdot}
\def\Si{\Sigma}
\def\lam{\lambda}
\def\bs{\backslash}
\def\l{^\ell}
\def\r{^r}
\def\dd{\ddot}
\def\dtm{\,\ddot{\tm}\,}
\def\dimp{\,\ddot{\imp}\,}
\def\de{\delta}
\def\cc{complemented category }
\def\ccs{complemented categories }
\def\cp{complemented pair }
\def\cps{complemented pairs }

\mathcaldef{C}
\mathcaldef{D}
\mathcaldef{U}
\mathcaldef{M}
\mathcaldef{S}
\mathcaldef{P}
\mathcaldef{B}
\mathcaldef{V}
\mathcaldef{L}
\mathcaldef{R}
\mathcaldef{A}
\mathcaldef{O}

\mathbfdef{Set}
\mathbfdef{Cat}
\mathbfdef{Pos}
\mathbfdef{Top}
\mathbfdef{nat}
\mathbfdef{ten}

\mathrmdef{id}
\mathrmdef{I}

\def\true{{\tt true}}
\def\false{{\tt false}}

\def\CatX{\Cat/X}
\def\PX{\Set^{X\op}}
\def\PY{\Set^{Y\op}}
\def\PXo{\Set^X}
\def\PYo{\Set^Y}
\def\PXX{\Set^{X\op\tm X}}
\def\PYY{\Set^{Y\op\tm Y}}
\def\VX{\V_0^{X\op}}
\def\VY{\V_0^{Y\op}}

\def\VXo{\V_0^X}
\def\VYo{\V_0^Y}
\def\co{\pi_0}

\def\hom{{\rm hom}}
\def\lim{{\rm lim}}
\def\colim{{\rm colim}}
\def\eend{{\rm end}}
\def\coend{{\rm coend}}
\def\Din{{\rm Din}}

\newtheorem{prop}{Proposition}
\newtheorem{corol}{Corollary}
\newtheorem{axiom}{Axiom}
\let\pf\proof
\let\epf\endproof
\def\eq{\begin{equation}}
\def\eeq{\end{equation}}

\author{Claudio Pisani}
\address{via Saluzzo 67,\\ 10125 Torino, Italy.}
\title{On indexed actions}
\copyrightyear{2011}
\eaddress{pisclau@yahoo.it}
\dedication{a mio padre}

\begin{document}

\maketitle

\begin{abstract}

We present some laws relating the $\Cat$-indexed categories 
of left, right and bi-actions:
by defining $(A\comp M)x = Mx^{Ax}$ one gets a biclosed monoidal action of $\PX$ on $(\PXo)\op$,
while $\B X$ and $\CatX$ act (partially) on their opposites by exponentials;
both the inclusions $(\B X,\B X)\to (\PX, \PXo) \to (\CatX,\CatX)$ preserve the (cartesian)
monoidal structures and the actions, and the same holds for substitutions along functors.
These strong morphisms of strong indexed monoidal actions have in fact a wider range of applications;
in particular, replacing $\Set$ with any (co)complete symmetric monoidal closed category $\V$,
we consider the pair of biclosed indexed monoidal actions $(\VX,\VXo;\comp\l_X,\comp\r_X;X\in\Cat)$ 
and its formal relationships with bi-actions and constant actions.

Some of the resulting laws also hold in a fragment of biclosed bicategory (with an object supporting a
symmetric monoidal category) and are taken, in the second part, as the basis for developing  
some abstract category theory. 
Finally, we add $\PXX$ to the picture and give a symmetrical version of the comprehension adjunction.
 
\end{abstract}








\section{Introduction}
\label{intro}

We are concerned with actions on a twofold level: both as the main object of study
and as the main tool to be used.
On the one hand, we are interested in actions of categories $X$ with respect to composition
(presheaves or more generally functors valued in a monoidal closed $\V$);
the relevant morphisms of $X$-actions are the usual ones: natural transformations. 
On the other hand, the actions of monoidal categories with respect to 
the (tensor) product also play a major role in our technical development; 
now the relevant morphisms involve the acting category as-well.

Consider the inclusion $i\l_X : \PX \to \CatX$ and $i\r_X : \PXo \to \CatX$ via discrete
(op)fibrations and let $\B X$ be their pullback, that is discrete bifibrations
with projections $j\l_X:\B X \to \PX$ and $j\r_X:\B X\to \PXo$ and $k_X := i\l_X j\l_X = i\r_X j\r_X$.
Exponentials (and products) in $\B X$ are computed as in $\CatX$:
$k_X(C^B) \iso (k_X C)^{k_X B}$. Otherwise stated, 
the pair $(k_X,k_X)$ (to be precise, $(k_X,k_X\op)$)
is an action morphism from the exponential action of $\B X$ on $(\B X)\op$ 
to the (partial) one of $\CatX$ on $(\CatX)\op$.
 If we define the action $\comp$ of $\PX$ on $(\PXo)\op$ pointwise by 
 \[
 (A\comp M) x := Mx^{Ax}  \qv (A\comp M)f := Mf\circ - \circ Af : Mx^{Ax} \to My^{Ay}
 \]
then $(k_X,k_X)$ factors through $(j\l_X,j\r_X)$ and $(i\l_X,i\r_X)$ which are also action morphisms:
 \[
 j\r_X(C^B) \iso j\l_X B \comp j\r_X C \qv i\r_X(A\comp M) \iso (i\r_X M)^{i\l_X A}
 \]
 Furthermore, any functor $f:X\to Y$ induces action morhisms at each of the three levels:
 \[
f\ast(C^B) \iso (f\ast C)^{f\ast B} \qv  f\r(A\comp M) \iso f\l A \comp f\r M \qv f\st(q^p) \iso (f\st q)^{f\st p}
 \]

A similar but not strictly action-like situation is the following.
Let $X$ be a topological space and denote by $\O X$, $\C X$ and $\B X$ the posets of 
open, closed and clopen parts respectively, with inclusions $j\l_X:\B X\to\O X$, $j\r_X:\B X\to\C X$,
$i\l_X:\O X \to \P X$ and $i\r_X:\C X \to \P X$.
The inclusion $k_X:\B X \to \P X$ is a Boolean algebra morphism, preserving in particular
the relative complement $A\imp B = \complement A \cup B$, which is the exponential action 
of $\B X$ on $(\B X)\op$ and of $\P X$ on $(\P X)\op$. 
The same formula gives an action $\comp : \O X \tm (\C X)\op \to (\C X)\op$, 
and the action morphism $(k_X,k_X)$ factors through the action morphisms $(j\l_X,j\r_X)$ and $(i\l_X,i\r_X)$. 
Again, any continuous map $f:X\to Y$ also induce action morphisms at the three levels.
The action $\comp : \O X \tm (\C X)\op \to (\C X)\op$ is biclosed:
the adjoint action $\O X \tm \C X \to \C X$ and the enrichment
$(\C X)\op \tm \C X \to \O X$ are given respectively by $\dm\r_X (i_X\l A\cap i_X\r M)$ 
(the \v closure" reflection in $\C X$ of the product in $\P X$) and by $\sq\l_X (i_X\r M \imp i_X\r N)$ 
(the \v interior" coreflection in $\O X$ of the exponential in $\P X$).

Since $A\comp M$ is, in disguise, simply the intersection of $A$ with the complement of $M$, 
the associated right adjoints are simply exponentials in $\O X$ in disguise.
While this is not the case for the set valued context, the formulas 
\[
\dm\r_X (i_X\l A\tm_X i_X\r M)  
\qv
\sq\l_X (i_X\r N)^{i_X\r M}
\] 
still give the associated right adjoints 
\[
\odot:\PX \tm \PXo \to \PXo
\qv
\tr : (\PXo)\op \tm \PXo \to \PX
\]
(where $i\l_X\adj \sq\l_X:\CatX\to\PX$ 
and $\dm\r_X \adj i\r_X:\PXo\to\CatX$).

In general we say that a biclosed monoidal action $\comp$ of $\V$ on $\M\op$ 
is a \v complemented category" and  a morphism of complemented categories 
is a pair of functors $f\l:\V\to\V'$ and $f\r:\M\to\M'$ which preserve the monoidal
structures and the actions and have right and left adjoint respectively. 
From $f\r(A\comp M) \iso f\l A \comp f\r M$, we then get by adjunction
$\ex\r_f(f\l A \odot M) \iso A\odot\ex\r_f M$, a sort of mixed Frobenius law.
For instance, for a topological space $X$ the morphism $(i\l_X,i\r_X):(\O X,\C X)\to (\P X,\P X)$ 
gives $\dm\r_X(i\l_X A \cap P) = A \odot \dm\r_X P = \dm\r_X(i\l_X A \cap i\r_X\dm\r_X P)$, 
which includes the fact that density is preserved on open parts.
Similarly, for a category $X$, the morphism $(i\l_X,i\r_X):(\PX,\PXo)\to (\CatX,\CatX)$ gives 
$\dm\r_X(i\l_X A \tm_X P) \iso A \odot \dm\r_X P \iso \dm\r_X(i\l_X A \tm_X i\r_X\dm\r_X P)$, which
is strictly related to the stability of final map with respect to discrete opfibrations (see~\cite{pis08}).  

As just sketched, left and right actions are united both by the (indexed) inclusion in $\CatX$ and 
by sharing the (indexed) subcategory $\B X$. 
(They are also united by the Isbell adjunction, so that they share $X$ as well, but this fact
does not seem to be strictly related to our present approach). 
While the first aspect can be useful for certain calculations, from an abstract point of view
it has the drawback that the categories $\CatX$ are not closed. 
Anyway, it is possible to capture the relevant formal laws of the second aspect by
taking in account the fact that the functors $\odot$ and $\tm$ on the one side and $\tr$ and $\comp$ 
on the other side collaps when one argument is restricted to $\B X$.
Thus we take the morphisms of indexed complemented categories
$(j\l_X,j\r_X):(\B X,\B X) \to (\PX,\PXo)$ as the basis for an abstraction in which
$\B 1$ has the role of \v internal truth values category": it turns out that all the categories 
and adjunctions involved in the definition are enriched in it.
The abstraction includes a sort of $\V$-relative category theory, for a symmetric monoidal
closed (co)complete category $\V$:
the left and the right actions of a category $X$ on $\V_0$ have monoidal structures 
induced pointwise by that of $\V$, and the action of each of them on the opposite
of the other one is induced pointwise by the internal hom of $\V$.

Further abstracting, we are naturally led to consider the concept of \v indexed pair" $(\L X,\R X;X\in\C)$ 
over a category $\C$ with a (not necessarly terminal) object $1\in\C$:
$\L 1$ and $\R 1$ are isomorphic and have a symmetric monoidal closed structure $\V$,
$\L X$ and $\R X$ have quantifications and are enriched, powered and copowered over $\V $;
furthermore, there are \v mixed tensor" bifunctors $*_X:\L X \tm \R X \to \V$ with enriched
\v absolute complement" adjoints $A*_X - \adj A\comp\l_X -:\V\to\R X$ 
and $-*_X M \adj M\comp\r_X -:\V\to\L X$,
and substitution functors preserve powers, copowers and complements:

\[     
\begin{array}{c}
\{V\otm\l_X A,B\}\l_X \\ \hline
\{V,\{A,B\}\l_X\} \\ \hline
\{A,[V,B]\l_X\}\l_X
\end{array}
\qv
\begin{array}{c}
\{V\otm\r_X M,N\}\r_X \\ \hline
\{V,\{M,N\}\r_X\} \\ \hline
\{M,[V,N]\r_X\}\r_X
\end{array}
\qv  
\begin{array}{c}
\{A *_X M,V\} \\ \hline
\{A,M\comp\r_X V\}\l_X \\ \hline
\{M,A\comp\l_X V\}\r_X 
\end{array}
\]
\[    
\begin{array}{c}
f\l(V\otm\l_Y A)  \\ \hline
V\otm\l_X f\l A  
\end{array}   
\qv
\begin{array}{c}
f\l[V,A]\l_Y  \\ \hline
[V,f\l A]\l_X   
\end{array}\qv
\begin{array}{c}
f\r (A\comp\l_Y V)  \\ \hline
f\l A \comp\l_X V  
\end{array}    
\]
\[     
\begin{array}{c}
f\r(V\otm\r_Y M)  \\ \hline
V \otm\r_X f\r M  
\end{array}   
\qv
\begin{array}{c}
f\r[V,M]\r_Y  \\ \hline
[V,f\r M]\r_X   
\end{array}\qv
\begin{array}{c}
f\l (M\comp\r_Y V)  \\ \hline
f\r M \comp\r_X V  
\end{array}    
\]
\[     
\begin{array}{c}
\{f\l A,B\}\l_X  \\ \hline
\{A,\fa\l_f B\}\l_Y
\end{array}   
\qv
\begin{array}{c}
\{A,f\l B\}\l_X  \\ \hline
\{\ex\l_f A,B\}\l_Y
\end{array}\qv
\begin{array}{c}
\{f\r M,N\}\r_X  \\ \hline
\{M,\fa\r_f N\}\r_Y
\end{array}   
\qv
\begin{array}{c}
\{M,f\r N\}\r_X  \\ \hline
\{\ex\r_f M,N\}\r_Y
\end{array}
\]

The same laws hold in a biclosed bicategory $\M$ (for instance, of $\V$-profunctors)
with a selected object $1$ which supports a {\em symmetric} monoidal category, 
by taking as $\C$ the \v maps" (right adjoint arrows) in $\B$ 
and posing $\L X := \B(X,1)$ and $\R X := \B(1,X)$.

In Section \ref{ip}, which can be red independently from the rest of the paper, 
we show that these axioms (along with some adequacy hypothesis) allow us to define weighted limits,
(pointwise) Kan extensions, fully faithful, dense and absolutely dense maps
and (if $1\in\C$ is actually terminal) conical limits and final maps
and to prove some of their familiar (and less familiar) properties with straightforward calculations. 

In the last section we came back to ordinary category theory, presenting
a generalization of the comprehension adjunction between categories over $X$
and presheaves on $X$ \cite{law70}.
The categories $X$ and $X\op$ are themselves united by their inclusions
in the groupoidal reflection $\ov X$, which induce the inclusions $j\l_X$ and $j\r_X$ 
of biactions in presheaves.
On the other hand, they can be also united by their product and the projections of $X\op\tm X$
induce the (dummy) inclusion of $\PX$ and $\PXo$ in $\PXX$; the discrete (op)fibrations
inclusions $i\l_X$ and $i\r_X$ factor through them and a \v diagonal" comprehension functor giving
the \v extension" $i_XH = \{x\in X|H(x,x)\}$ in $\CatX$ of the \v predicate" $H\in\PXX$.
In fact, the adjunctions $\dm\l_X\adj i\l_X :\PX \to \CatX$ and
$\dm\r_X\adj i\r_X :\PXo \to \CatX$ factor through $\dm_X\adj i_X :\PXX \to \CatX$,
where $\dm_X p \iso \ex_{p\op\tm p} \hom_P$, for $p:P\to X$.

The present paper is a development of previous work by the author
(see in particular~\cite{pis10} and references therein) but can be red independently.


\section{Complemented categories}
\label{comp}

Monoidal actions have been considered by several authors in different contexts.
In this section we show how they can be usefully seen has monoidal categories with a \v complement" functor
and present various instances of indexed monoidal actions.

\begin{definition}
A {\bf complemented category} $(\V, \M ) = (\V ,\otm, \I ; \M,\comp )$
is a symmetric monoidal closed category $(\V ,\otm, \I  )$ endowed with a {\bf complement} 
in a category $\M$, that is a biclosed monoidal action $\comp : \V \tm \M\op \to \M\op$
of $\V$ on $\M\op$. 
A {\bf morphism}  $(f\l,f\r):(\V,\M) \to (\V',\M')$ is a pair of functors which preserve
the monoidal structure and the action up to isomorphisms and which have a right and a left adjoint respectively.
\end{definition}
Thus, complemented categories and their morphisms are summarized (neglecting coherence and symmetry) 
by the following laws (natural isomorphisms):
\eq     \label{comp1}
\begin{array}{c}
\I\otm A \\ \hline
A
\end{array}
\qv
\begin{array}{c}
(A\otm B)\otm C \\ \hline
A \otm (B\otm C)
\end{array}
\qv  
\begin{array}{c}
\I\comp M \\ \hline
M
\end{array}
\qv
\begin{array}{c}
(A\otm B)\comp M \\ \hline
A \comp (B\comp M)
\end{array}
\eeq

\eq    \label{comp2}
\begin{array}{c}
\V(A\otm B,C) \\ \hline
\V(A,[B,C]) \\ \hline
\V(B,[A,C])
\end{array}
\qv
\begin{array}{c}
\M\op(A\comp M,N) \\ \hline
\M\op(M,A \odot N) \\ \hline
\V(A,N\tr M)
\end{array}
\iff  
\begin{array}{c}
\M(N,A\comp M) \\ \hline
\M(A \odot N,M) \\ \hline
\V(A,N\tr M)
\end{array}
\eeq

\eq     \label{comp3}
\begin{array}{c}
f\l \I \\ \hline
\I'
\end{array}
\qv
\begin{array}{c}
f\l(A\otm B) \\ \hline
f\l A\otm' f\l B 
\end{array}
\qv  
\begin{array}{c}
f\r(A\comp M) \\ \hline
f\l A \comp' f\r M
\end{array}
\eeq

\eq    \label{comp4}
\begin{array}{c}
\V'(f\l A,B) \\ \hline
\V (A, \fa_f\l B)
\end{array}
\qv
\begin{array}{c}
\M'(M,f\r N) \\ \hline
\M (\ex_f\r M,N)
\end{array}
\eeq

The laws~(\ref{comp1})  and~(\ref{comp3})  yield, by the adjunctions~(\ref{comp2}) and~(\ref{comp4}),
the equivalent ones 

\eq    \label{comp5}
\begin{array}{c}
[\I , A] \\ \hline
A
\end{array}
\qv
\begin{array}{c}
[A\otm B, C] \\ \hline
[A , [B, C]]
\end{array}
\eeq

\eq      \label{comp6}
\begin{array}{c}
\I \odot M \\ \hline
M
\end{array}
\qv
\begin{array}{c}
(A \otm B)\odot M \\ \hline
A \odot (B \odot M) \\ 
\end{array}
\qv
\begin{array}{c}
(A \odot M)\tr N \\ \hline
[A , M\tr N] \\ \hline
M\tr (A\comp N)
\end{array}
\eeq

\eq 
\begin{array}{c}   \label{comp7}
[A,\fa_f\l B] \\ \hline
\fa_f\l[f\l A,B]'
\end{array}
\qv
\begin{array}{c}
A \odot \ex_f\r M \\ \hline
\ex_f\r(f\l A \odot' M)  
\end{array}
\qv
\begin{array}{c}
\ex_f\r M\tr N \\ \hline
\fa_f\l(M\tr' f\r N)
\end{array}
\eeq

\begin{remarks}   \label{comp8}
Note in particular that 
\begin{enumerate}
\item
$\odot: \V\tm\M \to \M$ is itself a biclosed monoidal action of $\V$ on $\M$. 
\item
$\M$ is enriched (via $\tr$), powered and copowered over $\V$ (so as $\V$ itself);
indeed, the adjunctions relating $\comp$, $\odot$ and $\tr$ are enriched over $\V$ 
(so as those relating $\otm$ and $[-,-]$).
\item
Given a morphism $(f\l,f\r):(\V,\M) \to (\V',\M')$, $\V'$ and $\M'$ are also
enriched (via $\fa_f\l[-,-]'$ and $\fa_f\l(-\tr'-)$ respectively), powered and copowered over $\V$.
\item
The adjunctions $f\l\adj\fa_f\l$ and $\ex_f\r\adj f\r$ are also enriched over $\V$.
\end{enumerate}
\end{remarks}

\begin{examples}    \label{comp9h}
\begin{enumerate}
\item
Any symmetric monoidal closed category $\V$ gives rise to a \cc $(\V,\V)$, with $A\comp B := [A,B]$.
We say that $(\V,\M)$ is {\bf standard} if it is isomorphic to $(\V,\V)$.
Note that $(\V,\M)$ is standard iff there is an isomorphism $i:\M\to\V_0$ such that $i(-\comp -) \iso [-,i-]$.
A morphism of standard \ccs is essentially a strong morphism of monoidal categories preserving also the 
closed structure.

In particular, any Heyting algebra gives rise to a \cc where $A\comp B$
is the exponential $A\imp B$, that is the usual pseudocomplement of $A$ relative to $B$.
\item
Any symmetric monoidal closed category $\V$ gives rise to a \cc $(\V,\V\op)$, with $A\comp B := A\otm B$.
We say that $(\V,\M)$ is {\bf topological} if it is isomorphic to $(\V,\V\op)$.
Again, $(\V,\M)$ is topological iff there is an isomorphism $i:\M\to\V\op$ such that $i(-\comp -) \iso -\otm i-$.
A morphism of topological \ccs is essentially a strong morphism of monoidal categories.

In particular, any topological space $X$ gives rise to a topological \cc with $\V:=\O X$
(the Heyting algebra of open parts), $\M:=\C X \iso (\O X)\op$ (the poset of closed parts)  
and $A\comp M$ is given by the exponential $A\imp M$ in $\P X$,
that is the relative complement $\complement A \cup M$ in $\P X$. 
Any continuous map gives rise to a morphism of \ccs.
\item
If $\V$ is *-autonomous, then the standard $(\V,\V)$ is isomorphic to the topological $(\V,\V\op)$,
via the isomorphism $(-)\ast:\V\to\V\op$. (Conversely, if $(\V,\M)$ is both standard and topological, 
then $\V$ has a *-autonomous structure.)
\end{enumerate}
\end{examples}

Often complemented categories occur in a symmetrical fashion:

\begin{definition}
A {\bf complemented pair} $(\L, \R ) = (\L,\otm\l, \I\l ; \R, \otm\r,\I\r;\comp\l,\comp\r )$
consists of two symmetric monoidal closed categories, each one endowed with a complement 
in (the underlying category of) the other one; 
that is it consists of two complemented categories $(\L, \R_0)$ and $(\R, \L_0)$.
A {\bf morphism}  $(f\l,f\r):(\L,\R) \to (\L',\R')$ is a pair of functors such that both
$(f\l,f\r):(\L,\R_0) \to (\L',\R'_0)$ and $(f\r,f\l):(\R,\L_0) \to (\R',\L'_0)$
are morphisms of \ccs.
\end{definition}

Thus, half of the laws summarizing complemented pairs and their morphisms are  
\eq     \label{comp10}
\begin{array}{c}
\I\l\otm\l A \\ \hline
A
\end{array}
\qv
\begin{array}{c}
(A\otm\l B)\otm\l C \\ \hline
A \otm\l (B\otm\l C)
\end{array}
\qv  
\begin{array}{c}
\I\l\comp\l M \\ \hline
M
\end{array}
\qv
\begin{array}{c}
(A\otm\l B)\comp\l M \\ \hline
A \comp\l (B\comp\l M)
\end{array}
\eeq

\eq    \label{comp11}
\begin{array}{c}
\L(A\otm\l B,C) \\ \hline
\L(A,[B,C]\l) \\ \hline
\L(B,[A,C]\l)
\end{array}
\qv
\begin{array}{c}
\R\op(A\comp\l M,N) \\ \hline
\R\op(M,A \odot\l N) \\ \hline
\L(A,N\tr\l M)
\end{array}
\iff  
\begin{array}{c}
\R(N,A\comp\l M) \\ \hline
\R(A \odot\l N,M) \\ \hline
\L(A,N\tr\l M)
\end{array}
\eeq

\eq     \label{comp12}
\begin{array}{c}
f\l \I \l \\ \hline
{\I\l}'
\end{array}
\qv
\begin{array}{c}
f\l(A\otm\l B) \\ \hline
f\l A{\otm\l}' f\l B 
\end{array}
\qv  
\begin{array}{c}
f\r(A\comp\l M) \\ \hline
f\l A {\comp\l}' f\r M
\end{array}
\eeq

\eq    \label{comp13}
\begin{array}{c}
\L'(f\l A,B) \\ \hline
\L (A, \fa_f\l B) \\ \hline
\L (\ex_f\l A,B)
\end{array}
\qv
\begin{array}{c}
\R'(M,f\r N) \\ \hline
\R (M, \fa_f\r N) \\ \hline
\R (\ex_f\r M,N)
\end{array}
\eeq
and the other half is obtained by exchanging $\L$ and $\R$ and the superscripts $\l$ and $\r$.

\begin{definition}
A $\C-${\bf indexed \cc} (resp. {\bf pair})  
is a pseudo\-functor from $\C\op$ to the category of \ccs (resp. pairs):
\[
(\V X, \M X;X\in\C) = (\V_X,\otm_X, \I_X ; \M_X, \comp_X; X\in\C ) 
\]
\[
(\L X, \R X;X\in\C) = (\L_X,\otm_X\l, \I_X\l ; \R_X, \otm_X\r,\I_X\r;\comp_X\l,\comp_X\r; X\in\C ) 
\]
A morphism of $\C-$indexed \ccs (resp. pairs) is a family $(t\l_X,t\r_X;X\in\C)$ 
of morphisms of \ccs (resp. pairs) such that the obvious squares commute up to isomorphisms.
\end{definition}

The following proposition gives a standard way to construct $\Cat$-indexed \ccs (or pairs):

\begin{proposition}     \label{comp9b}
If $(\V,\M)$ is a complemented category with $\V_0$ complete and cocomplete
then, for any $X\in\Cat$, $(\V,\M)^X = (\V_0^X,\M^{X\op})$ has also a complemeted category 
structure and any functor $X\to Y$ gives rise to a morphism $(\V,\M)^Y \to (\V,\M)^X$.
If $(\V,\M)$ is topological, so it is also $(\V,\M)^X$. 
If $X$ is a groupoid and $(\V,\M)$ is standard, so it is also $(\V,\M)^X$.
\end{proposition}
\pf
The monoidal structure on $\V_0^X$ and the action on $(\M^{X\op})\op \iso (\M\op)^X$ are inherited
\v pointwise" by that of $(\V,\M)$.  
Thus $A\comp M$ is the diagonal $X\to X\tm X$ followed by $A\tm M$ 
and by the internal hom $[-,-]$ of $\V$.
The complement action is biclosed due to the (co)completeness of $\V$.
Any functor $f:X\to Y$ gives rise, via substitution, to a morphism of \ccs $(\V,\M)^Y \to (\V,\M)^X$,
due to the pointwise nature of the structural operations. 
The rest can be seen by a routine check.
\epf

We henceforth tacitly assume that the symmetric monoidal closed categories $\V$, $\L$ and $\R$ 
underlying the \ccs and the \cps are complete and cocomplete.

\begin{corol}    \label{comp9c}
Any \cc (resp. pair) gives rise to a $\Cat$-indexed \cc (resp. pair).
\epf
\end{corol}

\begin{examples}    \label{comp9e}
\begin{enumerate}
\item
A locally cartesian closed category $\C$ gives rise to the (standard) 
indexed \cp $(\C/X,\C/X;X\in\C)$.
Substitution along $f:X\to Y$ in $\C$ is the morphism of (cartesian and standard) \cps
$(f\st,f\st):(\C/Y,\C/Y)\to (\C/X,\C/X)$ given by pullback.
\item
A topos $\C$ gives rise to the (standard) indexed \cp $(\P X,\P X; X\in \C)$.
\item
By applying Corollary \ref{comp9c} to the standard \cp $(\V,\V)$,
we get the indexed \cp $(\VX,\VXo;X\in\Cat)$, with $A\comp\l_X M$ 
defined as in the proof of Proposition~\ref{comp9b} and $M\comp\r_X A$ symmetrically 
as the diagonal $X\op\to X\op\tm X\op$ followed by $M\op\tm A$ and by the internal hom $[-,-]$ of $\V$.
For any functor $f:X\to Y$, the morphism of \cps $(f\l,f\r):(\VY,\VYo) \to (\VX,\VXo)$
is obtained by substituting $f\op$ in $A:X\op\to\V$ and $f$ in $M:X\to\V$.
\item   
As an instance of example (3) above, for $\V = \Set$ we get the indexed
\cc $(\PX,\PXo;X\in\Cat)$, where $A\comp\l_X M : X \to \Set$ is defined \v pointwise" 
by $(A\comp\l_X M)x := \Set(Ax, Mx)$.
Any functor $f:X\to Y$ gives rise, via substitution, to a morphism $(\PY,\PYo) \to (\PX,\PXo)$.
If $X$ is a groupoid, then $(\PX,\PXo)$ is standard; indeed, the inverse functor $(-)\st:X\op\iso X$
induces $s:\PXo \iso\PX$, and $A\comp M$ gives the exponential in $\PX$ (modulo $s$):
\eq  \label{group}
s(A\comp M) \iso A\imp s M
\eeq 

The inclusions (via discrete fibrations and opfibrations)
$i_X\l:\PX\to \CatX$ and $i_X\r:\PXo\to \CatX$ form a morphism 
\[
(\PX,\PXo;X\in\Cat) \to (\CatX,\CatX;X\in\Cat)
\] 
of indexed \cps with a partial codomain (since $\CatX$ is not closed in general).
Indeed, $i_X\r(A\comp M)$ is the exponential $i_X\r M^{i_X\l A}$ in $\CatX$ and there are adjunctions
$\dm_X\l \adj i_X\l \adj \sq_X\l$ and $\dm_X\r\adj i_X\r \adj \sq_X\r$.

The pullback (intersection) $i_X\l \tm_X i_X\r$ gives the category $\B X$ of discrete bifibrations
with projections (inclusions) $j_X\l:\B X\to\PX$ and $j_X\r:\B X\to\PXo$ 
(as presheaves which act by bijections).
The indexed inclusion $(j_X\l,j_X\r;X\in\Cat)$ 
gives, by (\ref{group}), another instance of morphism of indexed complemented pairs:
\eq    \label{group2}
(\B X,\B X;X\in\Cat)\to(\PX,\PXo;X\in\Cat) 
\eeq
Composing  $(j_X\l,j_X\r)$ with $(i_X\l,i_X\r)$ one gets an inclusion $(k_X\l,k_X\r)$
\[
(\B X,\B X;X\in\Cat) \to (\CatX,\CatX;X\in\Cat) 
\]
of standard indexed \cps (with a partial codomain).
\item
A two-valued correspective of example (4) above is given by the
morphisms of (topological) indexed \cps
\[
(\B X,\B X;X\in\Pos) \to (\D X,\U X;X\in\Pos) \to (\P X,\P X;X\in\Pos)
\]
where $\D X$ (resp. $\U X$) are the down-closed (resp. up-closed) subsets of $X$,
$\P X$ are all subsets and $\B X \iso \P(\co X)$ are the up-down-closed subsets. 
The composite $(\B X,\B X;X\in\Pos) \to (\P X,\P X;X\in\Pos)$ 
is of course a morphism of indexed boolean algebras.
\item
Similarly, we have morphisms of (topological) indexed \cps
\[
(\B X,\B X;X\in\Top) \to (\O X,\C X;X\in\Top) \to (\P X,\P X;X\in\Top)
\]
where $\B X$ are the clopen subsets of $X$. 
\end{enumerate}
\end{examples}

The morphism of indexed complemented pairs (\ref{group2})
in fact can be extended to the more general context
of Example \ref{comp9e} (3), giving the inclusion
\[
(j\l_X,j\r_X;X\in\Cat) : (\B X,\B X;X\in\Cat) \to (\VX,\VXo;X\in\Cat)
\]
of those actions which act by invertible maps in $\V_0$.
Then it is easy to see that some of the operators on $(\VX,\VXo;X\in\Cat)$ 
collapse when applied to biactions.
We formalize this fact in the following

\begin{definition}
A morphism $(j_X\l,j_X\r;X\in\C) : (\B X,\B X;X\in\C) \to (\L X,\R X;X\in\C)$ with a standard domain
is said to endowe $(\L X,\R X;X\in\C)$ with {\bf biactions} if the following laws hold:
\eq   \label{comp20}
\begin{array}{c}
j\l_X V \otm\l_X A \\ \hline
j\r_X V \odot\r_X A  
\end{array}
\qv
\begin{array}{c}
M \comp\r_X j\l_X V  \\ \hline
M \tr\l_X j\r_X V  
\end{array}
\qv
\begin{array}{c}
 j\r_X V \otm\r_X M \\ \hline
j\l_X V \odot\l_X M 
\end{array}
\qv
\begin{array}{c}
A \comp\l_X j\r_X V  \\ \hline
A \tr\r_X j\l_X V  
\end{array}
\eeq
An object $X\in\C$ is {\bf groupoidal} if $j\l_X$ and $j\r_X$ are equivalences.
For brevity, we refer to an indexed \cp $(\L X,\R X;X\in\C)$ endowed with biactions and 
such that $\C$ has a groupoidal terminal object $1$ as a {\bf normal pair}. 
\end{definition}

Instances of normal pairs are thus $(\VX,\VXo;X\in\Cat)$ with the biactions inclusion, 
$(\D X,\U X;X\in\Pos)$ with the inclusion of up-down-closed sets
and $(\O X,\C X;X\in\Top)$ with the inclusion of clopen sets. 

\begin{remark}
Given an indexed \cp $(\L X,\R X;X\in\C)$ endowed with biactions and a groupoidal
object $G\in\C$, we get a normal pair $(\L X,\R X;X\in\C/G)$.
\end{remark}

\begin{remark}   \label{terminal}
If $(\L X,\R X;X\in\C)$ is an indexed \cp and $1\in\C$ is a terminal object,
by the remarks~\ref{comp8} it follows that: 
\begin{enumerate}
\item
All the categories $\L X$ are enriched over $\L 1$ 
and over $\R 1$ by $\fa_X\l[-,-]\l_X$ and $\fa_X\r(-\tr_X\r-)$ (where $X:X\to 1$).
They also have copowers, $(X\l -)\otm\l_X-$ and $(X\r -)\odot\r_X-$,
and powers $[X\l -,-]\l_X$ and $(X\r -)\comp\r_X-$.
Symmetrically, the $\R X$ are enriched, copowered and powered over $\R 1$  and over $\L 1$.
\item
The adjunctions $f\l\adj\fa_f\l$ and $\ex_f\l\adj f\l$ are also enriched respectively
over $\L 1$ and over $\R 1$ (and symmetrically for $\ex_f\r\adj f\r\adj\fa_f\r$).
\item
The adjunctions relating $\comp_X\l$, $\odot_X\l$ and $\tr_X\l$  
(so as those relating $\otm_X\l$ and $[-,-]_X\l$) are also enriched over $\L 1$ (and symmetrically).
\end{enumerate}
\end{remark}

If $(\L X,\R X;X\in\C)$ is a normal pair, then Remark \ref{terminal} 
and the laws~(\ref{comp20}) (or their adjoint ones) give that:

\begin{enumerate}
\item
All the categories $\L X$ and $\R X$ are enriched over $\B 1$ via
\[
\{ A,B \}\l_X :=
\begin{array}{c}
\j\l\fa\l_X [A,B]\l_X  \\ \hline
\j\r\fa\r_X(A \tr_X\r B)
\end{array}
\qv
\{ M,N \}\r_X :=
\begin{array}{c}
\j\r\fa\r_X [M,N]\r_X  \\ \hline
\j\l\fa\l_X(M \tr_X\l N)
\end{array}
\]
where $\j\l$ and $\j\r$ are adjoint to the equivalences $j\l = j_1\l$ and $j\r = j_1\r$.
They are also copowered and powered over $\B 1$:
\[
V\otm\l_X A :=
\begin{array}{c}
X\l j\l V \otm\l_X A \\ \hline
X\r j\r V \odot\r_X A
\end{array}
\qv
V\otm\r_X M :=
\begin{array}{c}
X\r j\r V \otm\r_X M \\ \hline
X\l j\l V \odot\l_X M 
\end{array}
\]
\[
[V,A]\l_X  :=
\begin{array}{c}
[X\l j\l V , A]\l_X  \\ \hline
X\r j\r V \comp\r_X A   
\end{array}
\qv
[V,M]\r_X  :=
\begin{array}{c}
[X\r j\r V , M]\r_X   \\ \hline
X\l j\l V \comp\l_X M   
\end{array}
\]
\item
There are {\bf absolute complement} and {\bf mixed tensor} adjoint functors
\[
A*_X M \adj_A A\comp_X\l V 
\qv
A*_X M \adj_M M\comp_X\r V 
\]
given by
\[
M \comp_X\r V :=
\begin{array}{c}
M \comp\r_X X\l j\l V  \\ \hline
M \tr\l_X X\r j\r V  
\end{array}
\qv
A \comp_X\l V :=
\begin{array}{c}
A \comp\l_X X\r j\r V  \\ \hline
A \tr\r_X X\l j\l V  
\end{array}
\]
\[
A*_X M :=
\begin{array}{c}
\j\l\ex\l_X ( M \odot\r_X A ) \\ \hline
\j\r\ex\r_X ( A \odot\l_X M )  
\end{array}
\]
\end{enumerate}

\subsection{Nine laws}
\label{nine}
Now we summarize the basic laws which relate the (left or right) \v actions" 
in a normal pair $(\L X,\R X;X\in\C)$ and the \v constant (bi)actions" in $\B 1$ 
(which can also be seen as the category of \v internal truth-values"),
with respect to a given map $f:X\to Y$ in $\C$.
We present them in three groups, such that the ones in the same group are each other
equivalent by adjunction (we omit the other nine obtained by left-right symmetry).

The first group says that substitution commutes with (or preserves) copowers, that the universal
quantification adjunction is enriched and that universal quantification commutes with powers:
\eq     \label{comp25}
\begin{array}{c}
f\l(V\otm\l_Y A)  \\ \hline
V\otm\l_X f\l A  
\end{array}   
\qv
\begin{array}{c}
\{ A,\fa\l_f B \}\l_Y \\ \hline
\{ f\l A, B \}\l_X
\end{array}    
\qv
\begin{array}{c}
[V,\fa\l_f B]\l_Y \\ \hline
\fa\l_f[V,B]\l_X 
\end{array}    
\eeq

The second group says that existential quantification commutes with copowers, that the existential
quantification adjunction is enriched and that substitution commutes with powers:
\eq       \label{comp26}
\begin{array}{c}
V \otm\l_Y (\ex\l_f A) \\ \hline
\ex\l_f (V\otm\l_X A)
\end{array}
\qv
\begin{array}{c}
\{ \ex\l_f A,B \}\l_Y \\ \hline
\{ A,f\l  B \}\l_X
\end{array}     
\qv
\begin{array}{c}
f\l[V,B]\l_Y  \\ \hline
[V,f\l B]\l_X   
\end{array}   
\eeq

The third group says that substitution can pass to the other argument inside a mixed tensor product
becoming an existential quantification, that substitution commutes with absolute complement
and that the absolute complement of an existentially quantified action is the same as 
the universal quantification of its absolute complement:
\eq       \label{comp27}
\begin{array}{c}
\ex\l_f A *_Y M \\ \hline
A *_X f\r  M
\end{array}
\qv
\begin{array}{c}
f\l (M\comp\r_Y V)  \\ \hline
f\r M \comp\r_X V  
\end{array}    
\qv
\begin{array}{c}
\ex\l_f A\comp\l_X V \\ \hline
\fa\r_f(A \comp\l_X V)  
\end{array}   
\eeq

\begin{remarks}  \label{normal}
\begin{enumerate}
\item
Besides the copowers - powers adjunction, also the mixed tensor - absolute complement
adjunction is enriched in $\B 1$. 
(In fact, since all the basic adjunctions defining a normal pair are enriched in $\B 1$, 
the same holds for the derived ones.) 
Explicitly, we have natural isomorphisms
\[
\begin{array}{c}
\{V\otm\l_X A,B\}\l_X \\ \hline
\{V,\{A,B\}\l_X\} \\ \hline
\{A,[V,B]\l_X\}\l_X
\end{array}
\qv
\begin{array}{c}
\{V\otm\r_X M,N\}\r_X \\ \hline
\{V,\{M,N\}\r_X\} \\ \hline
\{M,[V,N]\r_X\}\r_X
\end{array}
\qv  
\begin{array}{c}
\{A *_X M,V\} \\ \hline
\{A,M\comp\r_X V\}\l_X \\ \hline
\{M,A\comp\l_X V\}\r_X 
\end{array}
\]
\item
Most of the equations~(\ref{comp25}),~(\ref{comp26}) and~(\ref{comp27})
may be seen as expressing the fact that limits commute with limits, or that (co)limits can be defined 
in terms of limits (see the next section). 
On the other hand the first and the second groups may be seen as expressing the fact that
being a left (resp. right) adjoint is equivalent to preserving some kinds of colimits (resp. limits).
\item
The first group follows essentially from the fact that a morphism of \cps preserves the monoidal structures,
while the last two follow from the fact that it preserves the complement functors.
\end{enumerate}
\end{remarks}


\section{Some abstract category theory}
\label{ip}

In this section, which can be red independently from the rest of the paper, 
we develop some abstract category theory, resting on a few axioms which
hold true in a normal pair as well as in (a fragment of) a biclosed bicategory 
(for instance, that of $\V$-profunctors) with a suitable selected object 
(for instance the trivial $\V$-category; see Remark \ref{ipr1});
 these axioms allow us to define weighted limits,
(pointwise) Kan extensions, fully faithful, dense and absolutely dense maps
and (if $1\in\C$ is actually terminal) conical limits and final maps
and to prove (using also some adequacy hypothesis) 
some of their familiar (and less familiar) properties with straightforward calculations.

\subsection{Indexed pairs}

As we have seen in Section \ref{comp} (recall in particular Remark \ref{normal} (1)) 
any normal pair $(\L X,\R X;X\in\C)$ gives rise to an 
{\bf indexed pair} $(\L X,\R X;X\in\C)$ consisting of the following data and axioms:

\begin{enumerate}
\item
A category $\C$ and an object (not necessarly terminal) $1\in\C$.
\item
Two $\C$-indexed categories $\L X$ and $\R X$ with quantifications:
$\ex_f\l\adj f\l\adj\fa_f\l : \L X\to\L Y$, $\ex_f\r\adj f\r\adj\fa_f\r : \R X\to\R Y$,
for all $f:X\to Y$.
\item
$\L 1$ and $\R 1$ are isomorphic and have a symmetric monoidal closed structure 
$\V = (\V_0\iso\L1\iso\R1, \otm, I, [-,-])$.
\item
All the categories $\L X$ and $\R X$ are enriched ($\{A,B\}\l_X$ and $\{M,N\}\r_X$), powered 
($[V,A]\l_X$ and $[V,M]\r_X$) and copowered ($V\otm\l_X A$ and $V\otm\r_X M$) over $\V$.
\item
For any $X\in\C$ there is a \v mixed tensor" bifunctor $*_X:\L X \tm \R X \to \V$ with enriched
\v absolute complement" adjoints $A*_X - \adj A\comp\l_X -:\V\to\R X$ 
and $-*_X M \adj M\comp\r_X -:\V\to\L X$.
\item
All the operators collapse over $1\in\C$, becoming those of $\V$:
$\{V,W\}\l_1 \iso \{V,W\}\r_1 \iso [V,W]\l_1 \iso [V,W]\r_1 \iso [V,W] \iso V\comp\l_1 W \iso V\comp\r_1 W$
(and similarly for tensors).

In the following, we will thus use $\{-,-\}$ in place of $[-,-]$ for the internal hom of $\V$.
\item
Substitution functors preserve powers, copowers and complements.
\end{enumerate}

The laws (natural isomorphisms) which summarize an indexed pair are thus
(apart those concerning the associativity of $\otm$ and the quantification adjunctions):

\eq     \label{ip1}
\begin{array}{c}
\{V\otm\l_X A,B\}\l_X \\ \hline
\{V,\{A,B\}\l_X\} \\ \hline
\{A,[V,B]\l_X\}\l_X
\end{array}
\qv
\begin{array}{c}
\{V\otm\r_X M,N\}\r_X \\ \hline
\{V,\{M,N\}\r_X\} \\ \hline
\{M,[V,N]\r_X\}\r_X
\end{array}
\qv  
\begin{array}{c}
\{A *_X M,V\} \\ \hline
\{A,M\comp\r_X V\}\l_X \\ \hline
\{M,A\comp\l_X V\}\r_X 
\end{array}
\eeq

\eq     \label{ip2}
\begin{array}{c}
f\l(V\otm\l_Y A)  \\ \hline
V\otm\l_X f\l A  
\end{array}   
\qv
\begin{array}{c}
f\l[V,A]\l_Y  \\ \hline
[V,f\l A]\l_X   
\end{array}\qv
\begin{array}{c}
f\r (A\comp\l_Y V)  \\ \hline
f\l A \comp\l_X V  
\end{array}    
\eeq

\eq     \label{ip3}
\begin{array}{c}
f\r(V\otm\r_Y M)  \\ \hline
V \otm\r_X f\r M  
\end{array}   
\qv
\begin{array}{c}
f\r[V,M]\r_Y  \\ \hline
[V,f\r M]\r_X   
\end{array}\qv
\begin{array}{c}
f\l (M\comp\r_Y V)  \\ \hline
f\r M \comp\r_X V  
\end{array}    
\eeq

From the above laws, by adjunction, we obtain:

\eq     \label{ip4}
\begin{array}{c}
\{A,\fa\l_f B\}\l_Y \\ \hline
\{f\l A, B\}\l_X
\end{array}   
\qv
\begin{array}{c}
\{\ex\l_f A,B\}\l_Y \\ \hline
\{A,f\l  B\}\l_X
\end{array}    
\qv
\begin{array}{c}
\ex\l_f A *_Y M \\ \hline
A *_X f\r  M
\end{array}    
\eeq

\eq       \label{ip5}
\begin{array}{c}
[V,\fa\l_f A]\l_Y \\ \hline
\fa\l_f[V,A]\l_X 
\end{array}     
\qv
\begin{array}{c}
V \otm\l_Y (\ex\l_f A)  \\ \hline
\ex\l_f (V\otm\l_X A)
\end{array} 
\qv  
\begin{array}{c}
\ex\l_f A\comp\r_Y V \\ \hline
\fa\r_f(A \comp\r_X V)  
\end{array}
\eeq
(plus the symmetrical ones, obtained by exchanging $\l$ and $\r$).

\begin{remark}   \label{ipr1}
Given a biclosed bicategory $\M$ with a selected object $1$ such that $\V:=\M(1,1)$ is symmetric, 
we get an indexed pair by taking the maps (right adjoint arrows) in $\M$ as $\C$ 
and by posing $\L X := \M(X,1)$ and $\R X := \M(1,X)$.
Indeed, axioms (2) to (6) are a straightforward consequence of the closed structure of $\M$,
while axiom (7) is given by instances of associativity of composition in $\M$ (or of their adjoints),
when one of the arrows is a map.
Likewise, \v proarrow equipments" (\cite{wood}) give rise to indexed pairs.

Thus, what we do in this section can be also considered as an abstract approach
to category theory which uses only a fragment of the (abstract) profunctor category,
avoiding the full bicategorical machinary
(of which $\L X$ and $\R X$ and the various bifunctors are obviously a trace).
\end{remark}

\subsection{External weighted limits}

Recall that $A\in\L X$ (resp. $M\in\R X$) is to be thought of as a left (resp. right) action 
of $X$ on $\V$ (the trivial one, if $X=1$), that is as a ($\V$-)functor $X\op\to\V$ (resp. $X\to\V$). 
Thus $\{A,B\}\l\in\V$ (resp. $\{M,N\}\r\in\V$) is the \v internal truth value" of 
natural transformations $A\to B$ (resp. $M\to N$).
These \v $\V$-valued functors" can be \v composed" with the \v functors" $f:T\to X$ in $\C$,
giving $f\l A \in \L T$ and $f\r M \in \R T$;
in particular, if $T=1$ we get the value $x\l A\in\V$ (resp. $x\r M\in\V$) 
of $A$ (resp. $M$) at the \v point" $x:1\to X$. 
(Of course, the value at $x$ of $f\l A$ is the value at $fx$ of $A$.)

We can also \v compose" $A$ (or $M$) with some funtors $\V\to\V$;
namely, $[V,A]$ and $A\comp V$ can be seen as the substitution of $A$ in
the covariant and contravariant (enriched) functors represented by $V$. 
Indeed, posing $f = x:1\to X$ in (\ref{ip2}) 
\eq     
\begin{array}{c}
x\l[V,A]\l_X  \\ \hline
\{V,x\l A\}   
\end{array}\qv
\begin{array}{c}
x\r (A\comp\l_X V)  \\ \hline
\{x\l A, V\}  
\end{array}    
\eeq
we note that their value at $x$ is an internal hom of $\V$.

We now interpretate the indexed pair axioms in terms of \v external" (that is $\V$-valued) weighted limits. 
(The \v internal" ones are treated in Section \ref{internal}.)
We say (omitting the index $X$ when superfluous) that $\{A,B\}\l$ 
is the {\bf left limit} of $B$ weighted by $A$ (and similarly $\{M,N\}\r$ is a {\bf right limit}) 
and that $A*M$ is the {\bf colimit} of $M$ weighted by $A$ (or conversely).

The indexed pair laws
 \eq     \label{ipn1}
\begin{array}{c}
\{V,\{A,B\}\l\} \\ \hline
\{A,[V,B]\l\}\l
\end{array}
\qv
\begin{array}{c}
\{A * M,V\} \\ \hline
\{A,M\comp\r V\}\l 
\end{array}
\eeq
thus express limits and colimits in terms of limits or, more precisely, say that representables preserve 
limits and convert colimits into limits.

Similarly, the indexed pair laws
\eq      \label{ipn2}
\begin{array}{c}
[V,\fa\l_f A]\l_Y \\ \hline
\fa\l_f[V,A]\l_X 
\end{array}     
\qv
\begin{array}{c}
\ex\l_f A\comp\r_Y V \\ \hline
\fa\r_f(A \comp\r_X V)  
\end{array} 
\eeq
say that representables preserve external right (Kan) extensions and convert the left ones into right ones.

Using the (\ref{ip4}) for \v points" $x:1\to X$ we get the Yoneda and co-Yoneda isomorphisms:
\eq     \label{exy}
\begin{array}{c}
\{\ex\l_x \I,A\}\l_X \\ \hline
\{\I,x\l  A\}  \\ \hline
x\l A
\end{array}    
\qv
\begin{array}{c}
\ex\l_x \I *_X M \\ \hline
\I \otm x\r  M  \\ \hline
x\r M
\end{array}    
\eeq
and symmetrically, $\{\ex\r_x \I,M\}\r \iso x\r M$ and $A*\ex\r_x \I\iso x\l A$.

Thus, the \v images of points" are to be thought of as funtors (internally) represented by the point itself,
and a \v concrete representation" of the \v abstract category" $X\in\C$ is given by $\ov X\l$,
the full $\V$-enriched subcategory of $\L X$ generated by the $\ex\l_x\I$, for $x:1\to X$.
Since
\eq     
\begin{array}{c}
\{\ex\l_x \I,\ex\l_y \I\}\l \\ \hline
x\l  \ex\l_y \I \\ \hline
\ex\l_y \I*\ex\r_x \I \\ \hline
y\r  \ex\r_x \I \\ \hline
\{\ex\r_y \I,\ex\r_x \I\}\r
\end{array}    
\eeq
we see that $\ov X\l$ and $\ov X\r$ are dual.
(In fact, one should check composition.)

\subsection{Internal weighted limits}
\label{internal}

We define (internal) weighted limits by the internal correspective of the (\ref{ipn1}),
that is using the substitutions $f\l\ex\l_y \I$ and $f\r\ex\r_y \I$ in the internally 
representables functors. 
We say that $\{M,f\}:1\to Y$ is a {\bf limit} and that $A*f:1\to Y$ is a {\bf colimit}
of $f:X\to Y$, weighted by $M\in\R X$ and $A\in\L X$ respectively, if:
 \eq     
\begin{array}{c}
\{M,f\}\r\ex\r_y \I \\ \hline
\{M,f\r\ex\r_y \I\}\r
\end{array}
\qv
\begin{array}{c}
(A*f)\l\ex\l_y \I \\ \hline
\{A,f\l\ex\r_y \I\}\l
\end{array}
\eeq
(note that the naturality conditions in $y$ refer to the categories $\ov X\l$ and $\ov X\r$).
So, internal limits and colimits are defined in terms of external (left and right) limits.
\begin{remark}      \label{ref}
Now, the limit-colimit duality is a perfect symmetry, while it is not the case externally.
(One finds a similar situation in internal category theory.)
Note also that since  $\{M,f\}\r\ex\r_y \I = \ov Y\r(\ex\r_{\{M,f\}} \I, \ex\r_y \I)$,
the limit $\{M,f\}$, as a functor of its weight $M$, is a (partially defined) adjoint
of the restriction $\ov Y\r \to \R Y \to \R X$ of $f\r$;  
thus, it is also the (partially defined) reflection of $\ex\r_f M$ in $\ov Y\r$.
\end{remark}
The internal correspective of the laws
\eq     
\begin{array}{c}
f\r (A\comp\l V)  \\ \hline
f\l A \comp\l V  
\end{array}\qv
\begin{array}{c}
f\l[V,A]\l  \\ \hline
[V,f\l A]\l   
\end{array}    
\eeq
is simply the functoriality of substitution:
\eq     
\begin{array}{c}
f\r(g\r\ex_y\r \I)  \\ \hline
(gf)\r\ex\r_y \I   
\end{array}\qv
\begin{array}{c}
f\l(g\l\ex_y\l \I)  \\ \hline
(gf)\l \ex\l_y \I   
\end{array}    
\eeq
which (along with the adjunction $\ex\l_f \adj f\l$) allows us to get
\eq     
\begin{array}{c}
\{M,gf\}  \\ \hline
\{\ex\r_f M,g\}    
\end{array}\qv
\begin{array}{c}
A*gf  \\ \hline
\ex\l_f A*g    
\end{array}    
\eeq
that are the internal version of
\eq     
\begin{array}{c}
\{M,f\r  N\}\r   \\ \hline
\{\ex\r_f M,N\}\r 
\end{array}    
\qv
\begin{array}{c}
A * f\r  M  \\ \hline
\ex\l_f A * M 
\end{array}    
\eeq
Thus (since $\{\I,x\} = x = \I*x$) we also get the internal (co-)Yoneda isomorphisms:
\eq     
\begin{array}{c} 
\{\ex\r_x \I,f\}  \\ \hline
fx
\end{array}    
\qv
\begin{array}{c} 
\ex\l_x \I * f  \\ \hline
fx
\end{array}    
\eeq
The laws which pose quantifications on the right side
\eq     
\begin{array}{c}
\{f\r M,N\}\r   \\ \hline
\{M,\fa\r_f N\}\r 
\end{array}    
\qv
\begin{array}{c}
f\l A * M  \\ \hline
A * \ex\r_f M 
\end{array}    
\eeq
become internally the {\em definitions} of the {\bf Kan estensions} $\fa_f g$ and $\ex_f g$:
 \eq     \label{kan}
\begin{array}{c}
\{f\r M,g\}  \\ \hline
\{M,\fa_f g\}    
\end{array}\qv
\begin{array}{c}
f\l A*g  \\ \hline
A*\ex_f g    
\end{array}    
\eeq
They can be equivalently defined by external Kan extensions,
using the internal analogous of the (\ref{ipn2}):
\eq                \label{kan2}
\begin{array}{c}
(\fa_f g)\r\ex\r_x \I   \\ \hline
\fa_f\r(g\r\ex\r_x \I)  
\end{array}\qv
\begin{array}{c}
(\ex_f g)\l\ex\l_x \I   \\ \hline
\fa_f\l(g\l\ex\l_x \I)  
\end{array}    
\eeq
Indeed, there are natural isomorphisms
\[                
\begin{array}{c}
\{M,(\fa_f g)\r\ex\r_x \I\}   \\ \hline
\{M,\fa_f g\}\r\ex\r_x \I   
\end{array}\qv
\begin{array}{c}
\{M,\fa_f\r(g\r\ex\r_x \I)\}  \\ \hline
\{f\r M,g\r\ex\r_x \I \}  \\ \hline
\{f\r M,g\}\r\ex\r_x \I 
\end{array}    
\]




\subsection{The concrete representation of an indexed pair}
\label{concrete}

The concrete representation of $X\in\C$ as the $\V$-category $\ov X\l$ (or $\ov X\r$)
can be extended to \v abstract functors" $f:X\to Y$; indeed, since $\ex\l_f\ex\l_x\I = \ex\l_{fx}\I$,
we get a $\V$-functor $\ov f\l:\ov X\l \to \ov Y\l$ (or $\ov f\r:\ov X\r \to \ov Y\r$),
while for $A\in\L X$ we get a $\V$-funtor $\ov A\l:(\ov X\l)\op \to \V$ by restricting
$\{-,A\}\l$ to the \v representables" (and similarly for $M\in\R X$).
By Yoneda, each of them can be seen as the restriction to $\ov X\l$ or $\ov X\r$of (co)limits
as functors of their weight, in two ways (which should be equivalent via the duality $\ov X\l \iso (\ov X\r)\op$) :

\[ 
\xymatrix@1@C=3.5pc
{\ov A\l := \ov X\l \ar[r] & \L X \ar[r]^-{\{-,A\}\l} & \V} 
\qv
\xymatrix@1@C=3.5pc
{\ov A\r := \ov X\r \ar[r] & \R X \ar[r]^-{A*-} & \V} 
\]
\[ 
\xymatrix@1@C=3.5pc
{\ov M\r := \ov X\r \ar[r] & \R X \ar[r]^-{\{-,M\}\r} & \V} 
\qv
\xymatrix@1@C=3.5pc
{\ov M\l := \ov X\l \ar[r] & \L X \ar[r]^-{-*M} & \V} 
\]
\[ 
\xymatrix@1@C=3.5pc
{\ov f\r := \ov X\r \ar[r] & \R X \ar[r]^-{\{-,f\}} & \ov Y\r}\qv
\xymatrix@1@C=3.5pc
{\ov f\l := \ov X\l \ar[r] & \L X \ar[r]^-{-*f} & \ov Y\l} 
\]

The extent to which this representation is \v faithful" depends on the axioms
we discuss in the following section.

\subsection{Adequacy axioms}
\label{adq}

To develop abstract category theory in the frame of indexed pairs, we need
some \v reduction rules", that is some adequacy (or density) axioms which 
(along with Yoneda reduction itself) allow us to eliminate or introduce variables.

In the following, we have in mind as a model the category of small categories,
that is the indexed pair $(\PX, \PXo; X\in\Cat)$.
We do not consider here the question of the extent to which the axioms hold in $\V$-enriched contexts.

A (left or right) density condition for a functor $i$ with respect to a bifunctor $\cd $
has the form of a reduction rule for fractions:
\[
\begin{array}{c}
i x \cd y  \\ \hline
i x \cd y'   
\end{array}
\iff
\begin{array}{c}
y  \\ \hline
y'   
\end{array}
\qv
\begin{array}{c}
y \cd ix  \\ \hline
y' \cd ix   
\end{array}
\iff
\begin{array}{c}
y  \\ \hline
y'   
\end{array}
\]
Of course, we are adopting the convention that \v numerator" and \v denominator"
of a \v fraction" are to be intended as functors of the variables which appear in both of them, 
while the fraction itself indicates the existence of an isomorphism between them.
In fact, more generally, fractions are to be intended as morphisms (not necessarly isos)
whose relative direction depends on the variance of the arguments.

We will use the following very general principle of functorial calculus:

\begin{prop}    \label{fc}
Suppose that the bifunctors $f$ and $g$ are \v adjoint relatively to $h$ and $k$":
\[
\begin{array}{c}
h(f(x,y),z)  \\ \hline
k(x,g(y,z))   
\end{array}
\]
and that $i$ is left dense for $k$ and $j$ is right dense for $h$;
then $i$ is left dense for $f$ iff $j$ is right dense for $g$.
\end{prop} 
\pf
\[
\begin{array}{c}
f(ix,y)  \\ \hline
f(ix,y')
\end{array}
\iff
\begin{array}{c}
h(f(ix,y),jz)  \\ \hline
h(f(ix,y'),jz) 
\end{array}
\iff
\begin{array}{c}
k(ix,g(y,jz))  \\ \hline
k(ix,g(y',jz))   
\end{array}
\iff
\begin{array}{c}
g(y,jz)  \\ \hline
g(y',jz)   
\end{array}
\]
(We have supposed all functors covariant; otherwise, one has to do obvious changes.)

In the same way, when $h$ and $k$ are hom functors, $i$ and $j$ are identities and 
the intermediate category is $Y = 1+1$, one gets the usual \v mates" correspondence 
for ordinary adjunctions $f\adj g$ and $f'\adj g'$.
\epf

The (left and right) density of the identity with respect to the enriched hom of $\L X$ or $\R X$
is expressed by \v Yoneda reduction":
\eq     
\begin{array}{c}
A  \\ \hline
B 
\end{array}
\iff
\begin{array}{c}
\{C,A\}\l  \\ \hline
\{C,B\}\l  
\end{array} 
\iff
\begin{array}{c}
\{A,C\}\l  \\ \hline
\{B,C \}\l  
\end{array} 
\eeq

The first assumption says that the concrete representation of $A$ and $M$ 
as $\ov A\l$ and $\ov M\r$ are faithful, that is that the inclusions $\ov X\l \to \L X$ and $\ov X\r \to \R X$
of \v representables" in all \v presheaves" are dense:

\begin{axiom}       \label{ax1}
\eq          
\begin{array}{c}
A  \\ \hline
B 
\end{array}
\iff
\begin{array}{c}
\{\ex_x\l \I,A\}\l  \\ \hline
\{\ex_x\l \I,B\}\l  
\end{array}
\qv
\begin{array}{c}
M  \\ \hline
N 
\end{array}
\iff
\begin{array}{c}
\{\ex_x\r \I,M\}\r  \\ \hline
\{\ex_x\r \I,N\}\r  
\end{array}
\eeq
or equivalently:
\eq     
\begin{array}{c}
A  \\ \hline
B 
\end{array}
\iff
\begin{array}{c}
A * \ex_x\r \I  \\ \hline
B * \ex_x\r \I\end{array}
\qv
\begin{array}{c}
M  \\ \hline
N 
\end{array}
\iff
\begin{array}{c}
\ex_x\l\I * M  \\ \hline
\ex_x\l \I * N    
\end{array}
\eeq
\end{axiom}

As in the proof of Proposition \ref{fc}, Axiom \ref{ax1} implies the 
\v contraposition law" for absolute complement:

\eq             \label{ax1d}
\begin{array}{c}
A  \\ \hline
B 
\end{array}
\iff
\begin{array}{c}
B \comp V \\ \hline
A \comp V \end{array}
\qv
\begin{array}{c}
M  \\ \hline
N 
\end{array}
\iff
\begin{array}{c}
N \comp V  \\ \hline
M \comp V    
\end{array}
\eeq

Explicitly, we have:
\[
\begin{array}{c}
N \comp V  \\ \hline
M \comp V    
\end{array}
\iff
\begin{array}{c}
\{\ex\l_x\I , N \comp V \} \\ \hline
\{\ex\l_x\I , M \comp V \}
\end{array}
\iff
\begin{array}{c}
\{\ex\l_x\I * N , V \} \\ \hline
\{\ex\l_x\I * M , V \}
\end{array}
\iff
\begin{array}{c}
\ex\l_x\I * M  \\ \hline
\ex\l_x\I * N 
\end{array}
\iff
\begin{array}{c}
M  \\ \hline
N 
\end{array}
\]

Again by Proposition \ref{fc}, (\ref{ax1d}) is equivalent to

\eq              \label{ax1c}
\begin{array}{c}
A  \\ \hline
B 
\end{array}
\iff
\begin{array}{c}
A * M  \\ \hline
B * M  \end{array}
\qv
\begin{array}{c}
M  \\ \hline
N 
\end{array}
\iff
\begin{array}{c}
A * M  \\ \hline
A * N    
\end{array}
\eeq

We cannot assume a condition as Axiom \ref{ax1} for maps $f:X\to Y$ with respect to internal (co)limits, 
simply because $\C$ is not (yet) a 2-category.
Rather, we {\em define} such a structure on $\C$ by
\[
\C(X,Y)(f , g)  := \Cat(\ov X\l,\ov Y\l)(\ov f\l , \ov g\l ) \iso \Cat(\ov X\r,\ov Y\r)(\ov f\r , \ov g\r )
\]
that is (following our convention on fraction notation)
\eq     \label{ax2}
\begin{array}{c}
f  \\ \hline
g 
\end{array}
\iff
\begin{array}{c}
\ex_x\l\I * f  \\ \hline
\ex_x\l \I * g    
\end{array}
\iff
\begin{array}{c}
\{ \ex_x\r \I,f\}  \\ \hline
\{ \ex_x\r \I,g\} 
\end{array}
\eeq
Now, Axiom \ref{ax1} and (\ref{ax2}) imply
\eq     \label{ax2b}
\begin{array}{c}
f   \\ \hline
g
\end{array}
\iff
\begin{array}{c}
g\l \ex\l_y\I  \\ \hline
f\l \ex\l_y\I  
\end{array}
\iff
\begin{array}{c}
f\r \ex\r_y\I  \\ \hline
g\r \ex\r_y\I  
\end{array}
\eeq
Indeed (considering for instance the left hand side):
\[
\begin{array}{c}
g\l \ex\l_y\I  \\ \hline
f\l \ex\l_y\I  
\end{array}
\iff
\begin{array}{c}
\{\ex\l_x\I , f\l \ex\l_y\I \} \\ \hline
\{\ex\l_x\I , g\l \ex\l_y\I \}
\end{array}
\iff
\begin{array}{c}
(\ex\l_x\I * f )\l \ex\l_y\I \\ \hline
(\ex\l_x\I * g)\l \ex\l_y\I\end{array}
\iff
\begin{array}{c}
\ex\l_x\I * f  \\ \hline
\ex\l_x\I * g 
\end{array}
\iff
\begin{array}{c}
f  \\ \hline
g 
\end{array}
\]
where the third equivalence is Yoneda reduction in $\ov X\l$, 
since $x\l\ex\l_y\I \iso  \ov X\l(\ex\l_x\I,\ex\l_y\I)$ by (\ref{exy}).
Again by Proposition \ref{fc}, (\ref{ax2b}) is equivalent to 
\eq     
\begin{array}{c}
f  \\ \hline
g 
\end{array}
\iff
\begin{array}{c}
A * f  \\ \hline
A * g    
\end{array}
\iff
\begin{array}{c}
\{ M,f\}  \\ \hline
\{ M,g\} 
\end{array}
\eeq

As our second and last density condition, we assume:
\begin{axiom}
\eq     \label{ax3}
\begin{array}{c}
A  \\ \hline
B
\end{array}
\iff
\begin{array}{c}
A*_X f  \\ \hline
B*_X f
\end{array}
\qv
\begin{array}{c}
M  \\ \hline
N
\end{array}
\iff
\begin{array}{c}
\{ M,f \}_X  \\ \hline
\{ N, f \}_X
\end{array}
\eeq
(where naturality holds with respect to any category $\C(X,Y)$).
\end{axiom}

\subsection{Kan extensions}

We have defined Kan extensions in (\ref{kan}); 
from the properties of the 2-category structure of $\C$
it easily follows that they are defined up to isomorphisms:
if $f\iso f'$ and $g\iso g'$ then $\fa_f g \iso\fa_{f'}g'$ and $\ex_f g \iso\ex_{f'}g'$.
Furthermore, they are really extensions in $\C$:
\eq     
\begin{array}{c}
t \to \fa_f g  \\ \hline
t\r\ex\r_x\I \to (\fa_f g)\r\ex\r_x\I  \\ \hline
t\r\ex\r_x\I \to \fa_f\r (g\r\ex\r_x\I)  \\ \hline
f\r t\r\ex\r_x\I \to g\r\ex\r_x\I  \\ \hline
(tf)\r \ex\r_x\I \to g\r\ex\r_x\I  \\ \hline
tf \to g 
\end{array}
\qv
\begin{array}{c}
\ex_f g \to t  \\ \hline
t\l\ex\l_x\I \to (\ex_f g)\l\ex\l_x\I  \\ \hline
t\l\ex\l_x\I \to \fa_f\l (g\l\ex\l_x\I)  \\ \hline
f\l t\l\ex\l_x\I \to g\l\ex\l_x\I  \\ \hline
(tf)\l \ex\l_x\I \to g\l\ex\l_x\I  \\ \hline
g \to tf  
\end{array}
\eeq

\begin{remark}
Our Kan extensions are \v pointwise"; indeed they are preserved by rappresentables
(\ref{kan2}), and are given by the (co)limit formulas:
\eq     \label{kan3}
\begin{array}{c}
(\fa_f g)x  \\ \hline
\{\ex\r_x\I ,\fa_f g \}  \\ \hline
\{f\r\ex\r_x\I , g \}
\end{array}
\qv
\begin{array}{c}
(\ex_f g)x  \\ \hline
\ex\l_x\I *\ex_f g  \\ \hline
f\l\ex\l_x\I * g 
\end{array}
\eeq
\end{remark}

\subsection{Fully faithful maps}

Given a map $f:X\to Y$ in $\C$,
$\ov f\l:\ov X\l \to \ov Y\l$ is fully faithful iff
$\ov f\r:\ov X\r \to \ov Y\r$ is such and iff anyone of the following properties holds:
\eq
\begin{array}{c}
\{ \ex\l_x\I,\ex\l_y\I \}   \\ \hline
\{ \ex\l_{fx}\I,\ex\l_{fy}\I \} 
\end{array}
\qv
\begin{array}{c}
\{ \ex\r_x\I,\ex\r_y\I \}   \\ \hline
\{ \ex\r_{fx}\I,\ex\r_{fy}\I \} 
\end{array}
\qv
\begin{array}{c}
\ex\l_x\I * \ex\r_y\I   \\ \hline
\ex\l_{fx}\I * \ex\r_{fy}\I  
\end{array}
\eeq

\eq
\begin{array}{c}
\{A,B\}   \\ \hline
\{\ex\l_f A,\ex\l_f B\} 
\end{array}
\qv
\begin{array}{c}
\{M,N\}   \\ \hline
\{\ex\r_f M,\ex\r_f N\} 
\end{array}
\qv
\begin{array}{c}
\{M,g\}   \\ \hline
\{\ex\r_f M,\ex_f g\} 
\end{array}
\eeq

\eq
\begin{array}{c}
\{A,B\}   \\ \hline
\{\fa\l_f A,\fa\l_f B\} 
\end{array}
\qv
\begin{array}{c}
\{M,N\}   \\ \hline
\{\fa\r_f M,\fa\r_f N\} 
\end{array}
\qv
\begin{array}{c}
\{M,g\}   \\ \hline
\{\fa\r_f M,\fa_f g\} 
\end{array}
\eeq

\eq
\begin{array}{c}
A * M   \\ \hline
\ex\l_f A * \ex\r_f M
\end{array}
\qv
\begin{array}{c}
A * g   \\ \hline
\ex\l_f A * \ex_f g
\end{array}
\eeq

\eq
\begin{array}{c}
f\l\ex\l_f A   \\ \hline
A
\end{array}
\qv
\begin{array}{c}
f\r\ex\r_f M   \\ \hline
M
\end{array}
\qv
\begin{array}{c}
(\ex_f g)f   \\ \hline
g
\end{array}
\eeq

\eq
\begin{array}{c}
f\l\fa\l_f A   \\ \hline
A
\end{array}
\qv
\begin{array}{c}
f\r\fa \r_f M   \\ \hline
M
\end{array}
\qv
\begin{array}{c}
(\fa_f g)f   \\ \hline
g
\end{array}
\eeq

\eq
\begin{array}{c}
\ex\l_x \I * g   \\ \hline
\ex\l_{fx} \I * \ex_f g
\end{array}
\qv
\begin{array}{c}
f\l\ex\l_{fx}\I   \\ \hline
\ex\l_x\I
\end{array}
\qv
\begin{array}{c}
f\r\ex\r_{fx}\I   \\ \hline
\ex\r_x\I
\end{array}
\eeq

For the proof one uses the density rules (introduction and elimination of variables) 
of Section \ref{adq} and the adjunction-like laws, that is the introduction and elimination 
of quantifications rules summarized below for the reader convenience:
\eq     
\begin{array}{c}
\{f\l A, B\}\l   \\ \hline
\{A,\fa\l_f B\}\l 
\end{array}   
\qv
\begin{array}{c}
\{f\r M, N\}\r   \\ \hline
\{M,\fa\r_f N\}\r 
\end{array}    
\qv
\begin{array}{c}
\{f\r M, g\}  \\ \hline
\{M,\fa_f g\} 
\end{array}    
\eeq

\eq     
\begin{array}{c}
\{A,f\l  B\}\l  \\ \hline
\{\ex\l_f A,B\}\l 
\end{array}    
\qv
\begin{array}{c}
\{M,f\r  N\}\r  \\ \hline
\{\ex\r_f M,N\}\r
\end{array}    
\qv
\begin{array}{c}
\{M,gf\} \\ \hline
\{\ex\r_f M,g\}
\end{array}    
\eeq

\eq     
\begin{array}{c}
f\l A * M  \\ \hline
A * \ex\r_f M 
\end{array}  
\qv  
\begin{array}{c}
f\l A * g  \\ \hline
A * \ex_f g 
\end{array}  
\qv
\begin{array}{c}
A * f\r  M  \\ \hline
\ex\l_f A * M 
\end{array}  
\qv  
\begin{array}{c}
A * gf  \\ \hline
\ex\l_f A * g 
\end{array}  
\eeq

\pf
Let us prove some of the equivalences; the other ones can be proven with the
same technique and we leave them to the reader.
\[
\begin{array}{c}
\{ \ex\l_x\I,\ex\l_y\I \}   \\ \hline
\{ \ex\l_{fx}\I,\ex\l_{fy}\I \} 
\end{array}
\iff
\begin{array}{c}
\{ \ex\l_x\I,\ex\l_y\I \}   \\ \hline
\{ \ex\l_f\ex\l_x\I,\ex\l_f\ex\l_y\I \} 
\end{array}
\iff
\begin{array}{c}
\{ \ex\l_x\I,\ex\l_y\I \}   \\ \hline
\{ \ex\l_x\I,f\l\ex\l_f\ex\l_y\I \} 
\end{array}
\iff
\begin{array}{c}
\ex\l_y\I    \\ \hline
f\l\ex\l_f\ex\l_y\I 
\end{array}
\]

\[
\begin{array}{c}
\ex\l_x\I    \\ \hline
f\l\ex\l_f\ex\l_x\I 
\end{array}
\iff
\begin{array}{c}
\{ \ex\l_x\I,A \}   \\ \hline
\{ f\l\ex\l_f\ex\l_x\I,A \}
\end{array}
\iff
\begin{array}{c}
\{ \ex\l_x\I,A \}   \\ \hline
\{ \ex\l_x\I,f\l\fa\l_f A \} 
\end{array}
\iff
\begin{array}{c}
A    \\ \hline
f\l\fa\l_f A  
\end{array}
\]

\[
\begin{array}{c}
\ex\l_x\I    \\ \hline
f\l\ex\l_f\ex\l_x\I 
\end{array}
\iff
\begin{array}{c}
\ex\l_x\I * M   \\ \hline
f\l\ex\l_f\ex\l_x\I * M 
\end{array}
\iff
\begin{array}{c}
\ex\l_x\I * M   \\ \hline
\ex\l_x\I * f\r\ex\r_f M 
\end{array}
\iff
\begin{array}{c}
M    \\ \hline
f\r\ex\r_f M  
\end{array}
\]

\[
\begin{array}{c}
M    \\ \hline
f\r\ex\r_f M  
\end{array}
\iff
\begin{array}{c}
A * M   \\ \hline
A * f\r\ex\r_f M 
\end{array}
\iff
\begin{array}{c}
A * M   \\ \hline
f\l\ex\l_f A * M 
\end{array}
\iff
\begin{array}{c}
A   \\ \hline
f\l\ex\l_f A  
\end{array}
\]

\[
\begin{array}{c}
A   \\ \hline
f\l\ex\l_f A  
\end{array}
\iff
\begin{array}{c}
A * g   \\ \hline
f\l\ex\l_f A * g 
\end{array}
\iff
\begin{array}{c}
A * g   \\ \hline
A * (\ex_f g)f 
\end{array}
\iff
\begin{array}{c}
g   \\ \hline
(\ex_f g)f  
\end{array}
\]
and so on.
Note that in some cases one can chose different paths; for instance
\[
\begin{array}{c}
A    \\ \hline
f\l\fa\l_f A  
\end{array}
\iff
\begin{array}{c}
\{ B,A \}   \\ \hline
\{ B,f\l\fa\l_f A \} 
\end{array}
\iff
\begin{array}{c}
\{ B,A \}   \\ \hline
\{ f\l\ex\l_f B,A \}
\end{array}
\iff
\begin{array}{c}
B   \\ \hline
f\l\ex\l_f B
\end{array}
\]

\epf

\subsection{Absolutely dense maps}
\label{dense}

While fully faithful maps are those for which the unit of $\ex\l_f\adj f\l$ is an iso,
absolutely dense maps are those for which the counit of the same adjunction is an iso.
(In the context of a bicategory $\M$ of proarrows, these would become adjunctions {\em in} $\M$.) 
Absolutely dense, or \v connected", functors are treated (in an enriched context)
in \cite{dense} where one finds some of the characterizations below.

Given a map $f:X\to Y$ in $\C$,
the following properties are equivalent:
\eq
\begin{array}{c}
\{A,B\}   \\ \hline
\{f\l A,f\l B\} 
\end{array}
\qv
\begin{array}{c}
\{M,N\}   \\ \hline
\{f\r M,f\r N\} 
\end{array}
\qv
\begin{array}{c}
\{M,g\}   \\ \hline
\{f\r M,gf\} 
\end{array}
\qv
\begin{array}{c}
\{h,g\}^Y_Z   \\ \hline
\{hf,gf\} ^X_Z
\end{array}
\eeq

\eq
\begin{array}{c}
\ex\l_x\I * \ex\r_y\I   \\ \hline
f\l\ex\l_x\I * f\r\ex\r_y\I  
\end{array}
\qv
\begin{array}{c}
A * M   \\ \hline
f\l A * f\r M
\end{array}
\qv
\begin{array}{c}
A * g   \\ \hline
f\l A * gf
\end{array}
\eeq

\eq
\begin{array}{c}
\ex\l_f f\l A   \\ \hline
A
\end{array}
\qv
\begin{array}{c}
\ex\r_f f\r M   \\ \hline
M
\end{array}
\qv
\begin{array}{c}
\ex_f (gf)   \\ \hline
g
\end{array}
\eeq

\eq
\begin{array}{c}
\fa\l_f f\l A   \\ \hline
A
\end{array}
\qv
\begin{array}{c}
\fa \r_f f\r M   \\ \hline
M
\end{array}
\qv
\begin{array}{c}
\fa_f (gf)   \\ \hline
g
\end{array}
\eeq

\eq
\begin{array}{c}
\ex\l_f f\l\ex\l_x\I   \\ \hline
\ex\l_x\I
\end{array}
\qv
\begin{array}{c}
\ex\r_x f\r\ex\r_x\I   \\ \hline
\ex\r_x\I
\end{array}
\eeq
Proofs are very similar to those for fully faithful maps, and we leave them to the reader.

\begin{remark}   \label{surjective}
So as full faithfulness is a strong \v injectivity" property,
absolute density is a strong \v surjectivity" property.
We will treat presently other (weaker) surjectivity properties: left or right density and 
final or initial maps.

In fact, this is more than a vague analogy: in the indexed pair $(\P X,\P X; X\in\Set)$ 
we find again the usual concepts for mappings.
Note, by the way, that in that case $\L X = \R X$, $\V = {\bf 2} = \{\true, \false\}$,  
$\{P, Q\}$ and $P*Q$ are the truth values of $P\subseteq Q$ and of $P\cap Q \neq \emptyset$ 
respectively, and $P\comp\false$ is the usual complementary set.
The \v representables" become the sigletons and the weighted (co)limit $\{P,f\} = P*f$ 
exists iff $f$ is constant on $P$. 
\end{remark}

\subsection{Dense maps}

Given a map $f:X\to Y$ in $\C$,
the following proprieties are equivalent:
\eq    \label{dense2}
\begin{array}{c}
\{\ex\l_x\I, \ex\l_y\I\}   \\ \hline
\{f\l \ex\l_x\I,f\l \ex\l_y\I\} 
\end{array}
\qv
\begin{array}{c}
\fa\l_f f\l\ex\l_x\I   \\ \hline
\ex\l_x\I
\end{array}
\qv
\begin{array}{c}
\ex_f  f   \\ \hline
\id_Y
\end{array}
\qv
\begin{array}{c}
f\l\ex\l_x\I * f   \\ \hline
x 
\end{array}
\eeq
\pf
The equivalence between the first two conditions is readly obtained 
by eliminating $\ex_x\l\I$ on the left. Since the last one can be written
\( 
\begin{array}{c}
\ex\l_x\I * \ex_f f   \\ \hline
\ex\l_x\I * \id_Y
\end{array}
\)
we similarly obtain the equivalence between the last two.
The second one and the third one are equivalent because
\( 
\begin{array}{c}
 \fa_f\l f\l \ex\l_x\I  \\ \hline
(\ex_f f)\l \ex\l_x\I  
\end{array}
\).
\epf
A map satisfying these properties is said to be {\bf left dense}.
Indeed, by the last of (\ref{dense2}), $f:X\to Y$ is left dense iff any \v object"
of $Y$ is a colimit of $f$ \v canonically" weighted.
A map is {\bf right dense} if it satisfies the \v dual" (that is \v symmetrical") properties: 
\eq
\begin{array}{c}
\{\ex\r_x\I, \ex\r_y\I\}   \\ \hline
\{f\r \ex\r_x\I,f\r \ex\r_y\I\} 
\end{array}
\qv
\begin{array}{c}
\fa\r_f f\r\ex\r_x\I   \\ \hline
\ex\r_x\I
\end{array}
\qv
\begin{array}{c}
\fa_f  f   \\ \hline
\id_Y
\end{array}
\qv
\begin{array}{c}
\{ f\r\ex\r_x\I , f \}  \\ \hline
x 
\end{array}
\eeq
In the indexed pair $(\P X,\P X; X\in\Set)$ one finds again surjectivity.

\subsection{Limits preservation}

Given $f:X\to Y$, we say that $g:Y\to Z$ {\bf preserves the limit} $\{M,f\}$ if $\{M,gf\} \iso g\{M,f\}$
(and the same for colimits).
Similarly, if $h\fa_f g \iso \fa_f(hg)$, then we say that $h$ {\bf preserves the right Kan extension} $\fa_f g$
(and the same for the left ones).

We can now motivate the term \v absolute density" for the strong surjectivity notion
of Section \ref{dense}: the conditions
\eq
\begin{array}{c}
\ex_f (gf)   \\ \hline
g
\end{array}
\qv
\begin{array}{c}
\fa_f (gf)   \\ \hline
g
\end{array}
\eeq
imply (for $g=\id$) left and right density and show that left and right Kan extensions
\eq
\begin{array}{c}
\ex_f f   \\ \hline
\id_Y
\end{array}
\qv
\begin{array}{c}
\fa_f f   \\ \hline
\id_Y
\end{array}
\eeq
are {\bf absolute}, that is preserved by any map.
Furthermore, they are equivalent to
\eq
\begin{array}{c}
f\l\ex\l_x\I * gf   \\ \hline
g x
\end{array}
\qv
\begin{array}{c}
\{f\r\ex\r_x\I , gf\}   \\ \hline
g x
\end{array}
\eeq
that is to the fact that the density (co)limits
\eq
\begin{array}{c}
f\l\ex\l_x\I * f   \\ \hline
x
\end{array}
\qv
\begin{array}{c}
\{f\r\ex\r_x\I , f\}   \\ \hline
x
\end{array}
\eeq
are absolute.

\subsection{Conical limits}

Suppose now that $1\in\C$ is in fact terminal.
We can then define {\bf constant} \v functors" and \v presheaves"
as those that factor through $X:X\to 1$; thus
\[ 
yX:X\to Y \qv X\l V \in \L X \qv X\r V \in \R X 
\]
are the constant functors and presheaves whose values at
$x:1\to X$ are $yXx=y$, $x\l X\l V = V$ and $x\r X\r V = V$.

If we define $\I_X\l := X\l\I$ e $\I_X\r := X\r\I$, then
\eq   \label{con}
\begin{array}{c}
\ex\r_X M   \\ \hline
\I * \ex\r_X M  \\ \hline
X\l\I * M  \\ \hline
\I\l_X * M
\end{array}
\qv
\begin{array}{c}
\fa\l_X A   \\ \hline
\{\I , \fa\l_X A \} \\ \hline
\{ X\l\I , A \} \\ \hline
\{\I\l_X , A \}
\end{array}
\eeq
that is, external (co)limits weighted by the constant \v trivial actions" $\I\l_X$ 
give quantifications \v on all $X$".

A {\bf conical (co)limit} of $f:X\to Y$ is a (co)limit weighted by $\I_X$:
\[
\lim f := \{\I_X\r,f\}  \qv  \colim f := \I_X\l * f
\]
Conical (co)limits can be obtained as Kan extensions along $X\to 1$:
\eq
\begin{array}{c}
\ex_X f   \\ \hline
\colim f
\end{array}
\qv
\begin{array}{c}
\fa_X  f   \\ \hline
\lim f
\end{array}
\eeq
which is the internal correspective of (\ref{con}).
Indeed, considering for instance colimits and using (\ref{con}) itself:
\[
\begin{array}{c}
(\colim f)\l\ex\l_x\I  \\ \hline
\{\I\l_X, f\l\ex\l_x\I \} \\ \hline
\fa\l_X (f\l \ex\l_x\I)  \\ \hline
(\ex_X f)\l \ex\l_x\I
\end{array}
\]

Under the hypothesis that the canonical $\dm\l_X:\C/X \to \L X$ e $\dm\r_X:\C/X \to \R X$ 
have right adjoints $\dm\l_X\adj i\l_X$ and $\dm\r_X\adj i\r_X$
(the comprehension scheme of \cite{law70}) and that these are fully faithful, the weighted (co)limits
can be canonically reduced to conical (co)limits. 
Indeed, in that case any weight $A\in\L X$ is isomorphic to $\ex_t \I\l_T$ 
(for $t = i\l_X A : T\to X$) so that:
\[
\begin{array}{c}
A* f \\ \hline
\ex_t \I\l_T * f  \ \\ \hline
\I\l_T * ft  \\ \hline
\colim ft
\end{array}
\]
(Of course, if $\C = \Cat$ then $t = i\l_X A:T\to X$ is the discrete fibration associated to $A$.)

\subsection{Final maps}

Again in the hypothesis that $1\in\C$ is terminal, we have a further notion of \v surjectivity"; 
the following properties are equivalent for a map $f:X\to Y$ in $\C$:
\eq
\begin{array}{c}
\ex\l_f \I\l_X  \\ \hline
\I\l_Y  
\end{array}
\qv
\begin{array}{c}
\ex\r_X f\r M   \\ \hline
\ex\r_Y M 
\end{array}
\qv
\begin{array}{c}
\fa\l_X f\l A   \\ \hline
\fa\l_Y A 
\end{array}
\eeq

\eq
\begin{array}{c}
\ex\r_X f\r\ex\r_x \I  \\ \hline
\I 
\end{array}
\qv
\begin{array}{c}
\colim_X (gf)   \\ \hline
\colim_Y g 
\end{array}
\eeq
\pf 
First note that $\ex\l_X\ex\l_x\I = \ex\l_{Xx}\I = \I$. 
(If $\C=\Cat$, this corresponds to the fact that (the total of) the discrete fibration associated to a 
representable, that is a slice, is connected.)
Then, by (\ref{con}) we have:
\[
\begin{array}{c}
\ex\l_f \I\l_X  \\ \hline
\I\l_Y  
\end{array}
\iff
\begin{array}{c}
\ex\l_f \I\l_X * \ex\r_y\I\\ \hline
\I\l_Y * \ex\r_y\I
\end{array}
\iff
\begin{array}{c}
\I\l_X * f\r\ex\r_y\I\\ \hline
\ex\r_Y\ex\r_y\I 
\end{array}
\iff
\begin{array}{c}
\ex\r_X f\r\ex\r_y\I\\ \hline
\I
\end{array}
\]

\[
\begin{array}{c}
\ex\l_f \I\l_X  \\ \hline
\I\l_Y  
\end{array}
\iff
\begin{array}{c}
\ex\l_f \I\l_X * M \\ \hline
\I\l_Y * M
\end{array}
\iff
\begin{array}{c}
\I\l_X * f\r M \\ \hline
\I\l_Y * M\end{array}
\iff
\begin{array}{c}
\ex\r_X f\r M \\ \hline
\ex\r_Y M
\end{array}
\]

\[
\begin{array}{c}
\ex\l_f \I\l_X  \\ \hline
\I\l_Y  
\end{array}
\iff
\begin{array}{c}
\{\ex\l_f\I\l_X,A\}   \\ \hline
\{\I\l_Y,A\} 
\end{array}
\iff
\begin{array}{c}
\{\I\l_X,f\l A\}   \\ \hline
\{\I\l_Y,A\} 
\end{array}
\iff
\begin{array}{c}
\fa\l_X f\l A   \\ \hline
\fa\l_Y A 
\end{array}
\]

\[
\begin{array}{c}
\ex\l_f\I\l_X    \\ \hline
\I\l_Y  
\end{array}
\iff
\begin{array}{c}
\ex\l_f\I\l_X * g   \\ \hline
\I\l_Y * g 
\end{array}
\iff
\begin{array}{c}
\I\l_X * gf   \\ \hline
\I\l_Y * g 
\end{array}
\iff
\begin{array}{c}
\colim_X (gf)   \\ \hline
\colim_Y g 
\end{array}
\]

These maps are called {\bf final}, and symmetrically one defines {\bf initial} maps.

In the indexed pair $(\P X,\P X; X\in\Set)$ one finds again surjectivity.

\begin{remark}
Once we chose that the \v true" concrete representation of $X$ 
is, say, $\ov X\l$ (rather than $\ov X\r$), $A\in\L X$ is a \v presheaf" on $X$ 
(that is a \v contravariant functor" $X\to\V$) via the inclusion $\ov X\l \to \L X$, 
while $M\in\R X$ is a \v covariant functor" $X\to\V$.
So, reasoning \v on the left side", $f\l A$ is to be thought of as the substitution of $f\op$
(rather than $f$) in $A:X\op\to \V$. 
Thus the conditions
\[
\begin{array}{c}
\ex\r_X f\r M   \\ \hline
\ex\r_Y M 
\end{array}
\qv
\begin{array}{c}
\colim_X (gf)   \\ \hline
\colim_Y g 
\end{array}
\]
express the fact that precomposing with a final functor preserves colimits, while
\[
\begin{array}{c}
\fa\l_X f\l A   \\ \hline
\fa\l_Y A 
\end{array}
\]
expresses the fact that precomposing with $f\op$ (which is initial)
preserves limits.
\end{remark}

\begin{remark}
Since the absolute density of $f:X\to Y$ is equivalent to each one of the conditions:
\[
\begin{array}{c}
\ex\l_f f\l A    \\ \hline
A  
\end{array}
\qv
\begin{array}{c}
\ex\r_f f\r M    \\ \hline
M  
\end{array}
\]
and since $f\l\I\l_Y \iso \I\l_X$ and $f\r\I\r_Y \iso \I\r_X$, an absolutely dense map is both final and initial.
\end{remark}

\subsection{Adjunctible and adjoint maps}

A simple characterization of left adjoint functors $f:X\to Y$ in $\Cat$
is that, for all $y\in Y$, $f\l\ex\l_y\I \iso \ex\l_{f\st y}\I$, for a suitable $f\st y\in X$
(a reflection of $y$ along $f$). 
In the frame of indexed pairs, we call such a map {\bf left adjunctible},
rather than \v adjoint", since it does not implies the existence of a $g:Y\to X$
such that $g y \iso f\st y$.
Symmetrically one defines {\bf right adjunctible} maps.

\begin{remarks}
\begin{enumerate}
\item
In the indexed pair $(\P X,\P X; X\in\Set)$ one finds bijections
(see Remark \ref{surjective}).
\item
If $1\in\C$ is terminal, a left (resp. right) adjunctible map is initial (resp. final).
Indeed $\ex\l_X f\l\ex\l_y\I \iso \ex\l_X\ex\l_{f\st y}\I \iso \I$.
\end{enumerate}
\end{remarks}

Right (resp. left) adjunctible maps preserve limits (resp. colimits)
and right (resp. left) Kan extensions:
\[
\begin{array}{c}
(f\{M,g\})\r\ex\r_y\I    \\ \hline
\{M,g\}\r f\r\ex\r_y\I    \\ \hline
\{M,g\}\r \ex\r_{f\st y}\I    \\ \hline
\{M,g\r \ex\r_{f\st y}\I\}    \\ \hline
\{M,g\r f\r\ex\r_y\I\}    \\ \hline
\{M,(fg)\r\ex\r_y\I\}    \\ \hline
\{M,fg\}\r\ex\r_y\I    
\end{array}
\qv
\begin{array}{c}
(f\fa_h g)\r\ex\r_y\I    \\ \hline
(\fa_h g)\r f\r\ex\r_y\I    \\ \hline
(\fa_h g)\r \ex\r_{f\st y}\I    \\ \hline
\fa_h\r (g\r \ex\r_{f\st y}\I )    \\ \hline
\fa_h\r (g\r f\r\ex\r_y\I )    \\ \hline
\fa_h\r((fg)\r\ex\r_y\I )   \\ \hline
(\fa_h fg)\r\ex\r_y\I    
\end{array}
\]

It seems natural to define (left and right) {\bf adjoint} maps by their classical characterization 
in terms of Kan extensions: $g\adj f$ iff
\eq
\begin{array}{c}
tg    \\ \hline
\fa_f t 
\end{array}
\iff
\begin{array}{c}
\{M,tg\}    \\ \hline
\{M,\fa_f t\} 
\end{array}
\iff
\begin{array}{c}
\{\ex\r_g M,t\}    \\ \hline
\{f\r M, t\} 
\end{array}
\iff
\begin{array}{c}
\ex\r_g M    \\ \hline
f\r M 
\end{array}
\eeq
which, since $\ex\r_g \adj g\r$, is equivalent to $f\r\adj g\r$.
Symmetrically, 
\(
\begin{array}{c}
tf    \\ \hline
\ex_g t 
\end{array}
\)
is equivalent to $g\l \adj f\l$. 


Of course, a left (resp. right) adjoint map is left (resp. right) adjunctible:
$f\r\ex\r_y\I \iso \ex\r_{g y}\I$.


%



\section{The symmetrical comprehension adjunction}
\label{coend}

In this section we consider the $\Cat$-indexed category $\PXX$,
with substitution along $f:X\to Y$ given by $\dd{f}H(x,x') := H(fx,fx')$,
and its relationships with $\PX$, $\PXo$ and $\CatX$.
While we do not propose here any abstraction or generalization, some of the results presented
may suggest steps in that direction.

\subsection{The indexed category of endoprofunctors}

Note that $(\PXX;X\in\Cat)$ can be seen as $(\PXo;X\in\Cat)$ restricted to
the categories of the form $X\op\tm X$ and to the functors of the form $f\op\tm f$,
so that $\dd f$ becomes $(f\op\tm f)\r$ and has adjoints $\dd\ex_f\adj \dd f \adj \dd\fa_f$
given by $\dd\ex_f := \ex\r_{f\op\tm f}$ and $\dd\fa_f := \fa\r_{f\op\tm f}$.
Since $X\op\tm X$ is canonically self-dual, any one of its right actions $H$ correspond to
a left action $H'$ given by $H'(x,y) := H(y,x)$. 
Furthermore, $(f\tm f\op)\l H' = ((f\op\tm f)\r H)'$.
The projections of $X\op\tm X$ induce the \v dummy inclusion" indexed functors 
$\de\l:\PX \to \PXX$, $\de\r:\PXo \to \PXX$.

Given $A:X\op\to\Set$ and $M:X\to\Set$, we define $A\dtm M$ 
by $(A\dtm M)(x,y) = Ax \tm My$, that is as the composite 
\[ \xymatrix@1@C=3.5pc{X\op\times X \ar[r]^-{A\times M} & \Set\times\Set \ar[r]^-\times & \Set} \]
This \v operation" $\dtm:\PX \tm \PXo \to \PXX$ is indexed, in the sense that $f\l A \dtm f\r M \iso \dd{f}(A\dtm M)$.
(Note that $A\dtm M$ is the product $\de\l A \tm \de\r M$ in $\PXX$.)
Another operation $\dimp\r : \PXo\tm\PXo\to\PXX$ is obtained by posing $(M\dimp\r N)(x,y) = [Mx, Ny]$:
\[ \xymatrix@1@C=3.5pc{X\op\times X \ar[r]^-{M\times N} & \Set\op\times\Set \ar[r]^-{[-,-]} & \Set} \]
(and similarly one defines $\dimp\l : \PX\tm\PX\to\PXX$).
Note that $\I_X \dtm M \iso \de\r_X M \iso \I_X \dimp\r M$, $A\dtm \I_X \iso \de\l_X A \iso \I_X\dimp\l A$,
where $\I_X$ is the terminal presheaf on $X$.

Now, given $H\in\PXX$, let us define a category $i_XH$ over $X$ as follows:
\begin{itemize}
\item
the objects over $x\in X$ are the elements of $H(x,x)$;
\item
given $\lam:x\to y$ in $X$, there is at most one arrow from $a\in H(x,x)$ to $b\in H(y,y)$ over $\lam$,
and this is the case iff $H(x,\lam)a = H(\lam,y)b \in H(x,y)$. 
\end{itemize}
Then one easily verifies that (see also \cite{pis07} and the references therein): 
\begin{enumerate}
\item
This constructions is the object map of an indexed functor
$i_X:\PXX\to\CatX$.
\item
The \v category of elements" functors $i_X\l:\PX\to\CatX$ and $i_X\l:\PXo\to\CatX$ factor through it:
$i_X\l \iso i_X\de\l$, $i_X\r \iso i_X\de\r$.
\item
$i_X(A\dtm M)$ is the product $i\l_X A\tm i\r_X M$ in $\CatX$. 
In particular, for objects $x,y: 1\to X$, $i_X(\ex\l_y\I\dtm \ex\r_x\I) \iso x\bs X\tm X/y$ is the
\v interval" category $[x,y]$ (over $X$) with objects $x\to z\to y$.
\item
$i_X(M\dimp\r N)$ is the exponential $(i\r_X N)^{i\r_X M}$ in $\CatX$
(and similarly for $i_X(A\dimp\l B)$).
\item
The inclusion $i_X:\PXX\to\CatX$ is not full.
Indeed, 
\[
\CatX(i_X H,i_X K) \iso\Din^*_X(H,K)
\]
 that is one gets the strong dinatural transformations, also known as \v Barr dinatural". 
Recall also that 
\[  
\Din_X(A\dtm M,K) \iso \Din^*_X(A\dtm M,K) 
\]
that is, the dinatural transformations with domain $A\dtm M$ are also strongly dinatural, for any $K$.
\end{enumerate}

Then, since the end and the (strong) coend of $H$ are representations of the functors
$\Din_X(\Delta_X S,H)$, $\Din_X^*(H,\Delta_X S)$ and $\Din_X(H,\Delta_X S)$
respectively, we get:
\begin{proposition}
The set of sections of $i_X H$
\[
\CatX(\id_X,i_X H) \iso \fa\l_X\sq\l_X(i_X H) \iso \fa\r_X\sq\r_X(i_X H)
\]
gives the end of $H:X\op\tm X \to \Set$.
The components of (the total of) $i_X H$
\[
\dm\Si_X(i_X H) \iso \ex\l_X \dm\l_X(i_X H) \iso \ex\r_X \dm\r_X(i_X H)
\]
gives the strong coend of $H$. The coend of $A\dtm M$ coincides with its strong coend
and with the mixed tensor product of Section \ref{comp}:
\[
A * M \iso \dm\Si_X (i\l_X A \tm i\r_X M) \iso \ex\l_X\dm\l_X (i\l_X A \tm i\r_X M) \iso \ex\r_X\dm\r_X(i\l_X A \tm i\r_X M)
\]
\epf
\end{proposition}

\begin{remarks}   \label{r0}
\begin{enumerate}
\item
It is natural to see $i_X H$ as a sort of \v diagonal" extension of $H$ and to denote it
by $\{x\in X|H(x,x)\}$ (see~\cite{law70} and Remark~\ref{tv}).
Then the usual end {\em notation} $\int_{x\in X}H(x,x)$ can be replaced (for set-valued functors)
by the end {\em formula} $\int_X\{x\in X|H(x,x)\}$ of the above proposition, where $\int_X$ 
is the sections functor, right adjoint to $\Set\to\CatX$.
\item
As a corollary, one gets the naturality formula for set-valued functors:
\[
\begin{array}{c}
\eend(M\dimp\r N)  \\ \hline
\CatX(\id_X,i_X(M\dimp\r N))  \\ \hline 
\CatX(\id_X,(i\r_X N)^{i\r_X M})  \\ \hline
\CatX(i\r_X M,i\r_X N)  \\ \hline
\PXo(M,N)
\end{array}
\]
\item
While strong dinaturality has the advantage over dinaturality of arising naturally as a full 
subcategory of $\CatX$ (so that, for instance, strong dinatural transformations always compose)
there are other facts that conversely seem to indicate a prevalent role for dinaturality.
For instance $\Din_X(H,K)$ (and not $\Din^*_X(H,K)$) can be expressed as an end.
Furthermore, the formula (\ref{dy}) below has a correspective for coends:
the mixed tensor product $H' * \hom_X$ ($\iso \hom_{X\op} * H$) 
gives the (not strong) coend of $H$.
\end{enumerate}
\end{remarks}

The following (rephrased) is referred to, in~\cite{mac65}, as \v diagonal Yoneda":
\eq      \label{dy}
\CatX(\id_X,i_X H) \iso \PXX(\hom_X, H)
\eeq
and says that there are two ways to express the end of an endoprofunctor.
In fact, it is also the key fact to prove the
\begin{prop}  \label{p1}
$i_X:\PXX\to\CatX$ has a left adjoint $\dm_X\adj i_X$, which takes $p:P\to X$ to 
$\dd\ex_p\hom_P$.
\end{prop}
\pf
\[
\begin{array}{c}
\PXX(\dd{\ex}_p\hom_P,H) \\ \hline
\Set^{P\op\tm P}(\hom_P,\dd{p}H) \\ \hline
\Cat/P(\id_P,i_P(\dd{p}H)) \\ \hline
\Cat/P(\id_P,p\st(i_X H)) \\ \hline
\CatX(p,i_X H)\end{array}
\]
\epf
In particular, $\hom_X$ is the reflection of the terminal $\id_X$ of $\CatX$.

The value of $\dm_X p$ at $\la x,y\ra$ can be expressed in various ways:
\eq      \label{r1}
\begin{array}{c}
(\dm_X p)(x,y) \\ \hline
p\l\ex\l_y\I * p\r\ex\r_x\I \\ \hline
\coend_P (p\l\ex\l_y\I\dtm p\r\ex\r_x\I) \\ \hline
\coend_P ^*(p\l\ex\l_y\I\dtm p\r\ex\r_x\I) \\ \hline
\coend_P^*\dd{p}(\ex\l_y\I\dtm \ex\r_x\I) \\ \hline
\dm \Si_P i_P\dd{p}(\ex\l_y\I\dtm \ex\r_x\I) \\ \hline
\dm \Si_P p\st i_X(\ex\l_y\I\dtm \ex\r_x\I) \\ \hline
\dm \Si_P p\st [x,y] \\ \hline
\dm \Si_X(p \tm [x,y]) 
\end{array}
\eeq
Indeed, using the usual formula for left Kan extensions (see also Section \ref{ip})
\[ (\ex\r_f M)x \iso f\l\ex\l_x\I * M \]
and observing that, for $\la x,y\ra : 1 \to X\op\tm X$,  $\ex\l_{\la x,y\ra}\I \iso \ex\r_x\I\dtm \ex\l_y\I$, we get 
\[ 
(\ex\r_{p\op\tm p} \hom_P)(x,y) \iso (p\op\tm p)\l\ex\l_{\la x,y\ra}\I * \hom_P \iso (p\r\ex\r_x\I\dtm p\l\ex\l_y\I) * \hom_P
\]
which is the third row above (see Remark \ref{r0} (3)). The other equivalences are immediate.

\begin{remark}
Since $i\r_X \iso i_X\de\r_X$ and $\de\r_X \iso \pi_2\r$ for the projection $\pi_2:X\op\tm X\to X$,
we also have $\dm\r_X \iso \ex\r_{\pi_2}\dm_X$, that is
$\ex\r_p\I_P \iso \ex\r_{\pi_2}\ex\r_{p\op\tm p}\hom_P \iso \ex\r_p\ex\r_{\pi_2}\hom_P$, for any $p\in\CatX$. 
Indeed, $\ex\r_{\pi_2}\hom_P \iso \I_P$, as one can easily verify directly:
\[ 
(\ex\r_{\pi_2} \hom_P)x \iso \pi_2\l\ex\l_x\I * \hom_P \iso \coend(\de\l\ex\l_x\I) \iso 
\ex\l_X\dm\l_X i_X(\de\l\ex\l_x\I) \iso \ex\l_X \dm\l_X i\l_X \ex\l_x\I \iso \I
\]
\end{remark}

\begin{remark}  \label{tv}
Classically, the third row of (\ref{r0}) is written $\int^{a\in P} X(x,pa)\tm X(pa,y)$.
As argued elsewhere, the last row of (\ref{r0}) can be seen as the set-valued version
of the predicate $p$ \v meets" $[x,y]$. 
In fact, in two valued contexts we have simplified forms of the adjunction $\dm_X\adj i_X$.  
For instance, let $X$ be a poset, $2^{X\op\tm X}$ the poset of binary relations on $X$
compatible with the order, and $\P X$ the poset of all parts of $X$. 
Then we have $\dm_X\adj i_X:2^{X\op\tm X}\to\P X$, where $i_X H = \{x\in X | H(x,x)\}$ and
$x(\dm_X P)y \iff \ex a\in P(x\leq a \leq y) \iff P\cap[x,y] \neq \emptyset$,
where $[x,y]= \{ z\in X | x\leq z \leq y \}$; indeed, as before,  $[x,y] = i_X H_{xy}$ where
$H_{xy}$ is the product of representables of opposite variance: $z H_{xy} w \iff x\leq w \,\&\, z\leq y$.
\end{remark}

\subsection{Absolutely dense and fully faithful functor}

We conclude by adding further characterizations of absolutely dense 
and of fully faithful functors (see Section \ref{ip}).

\begin{prop}    \label{ad}
Each of the following is equivalent to the absolute density of the functor $f:X\to Y$:
\begin{enumerate}
\item
$\ex_{f\op\tm f}\hom_X \iso \hom_Y$;
\item
$\eend_Y H \iso \eend_X\dd f H$, naturally in $H$; 
\item
$\coend_Y H \iso \coend_X\dd f H$, naturally in $H$; 
\item
$\dm_Z g \iso \dm_Z(gf)$, for any $g:Y\to Z$; 
\item
$\dm_Y f \iso \hom_Y$;
\item
$g$ and $g(f\op\tm f)$ have the same (co)end, for any $g:Y\op\tm Y\to Z$. 
\end{enumerate}
\end{prop}
\pf
The first condition is equivalent to absolute density by (\ref{r1})
(see Section \ref{dense}) and is equivalent to the second one by
\[
\begin{array}{c}
\eend_X \dd f H \\ \hline
\PXX(\hom_X, \dd f H) \\ \hline
\PYY(\dd\ex_f\hom_X, H) \\ \hline
\PYY(\hom_Y, H) \\ \hline
\eend_Y H
\end{array}
\qv
\begin{array}{c}
\PYY(\dd\ex_f\hom_X, H) \\ \hline
\PXX(\hom_X, \dd f H) \\ \hline
\eend_X \dd f H \\ \hline
\eend_Y H  \\ \hline
\PYY(\hom_Y, H) 
\end{array}
\]
Similarly one proves the equivalence between the first one and the third one:
\[
\begin{array}{c}
\coend_X \dd f H \\ \hline
\hom_{X\op} * \dd f H \\ \hline
\dd\ex_{f\op}\hom_{X\op} * H \\ \hline
\hom_{Y\op} * H \\ \hline
\coend_Y H
\end{array}
\qv
\begin{array}{c}
\dd\ex_{f\op}\hom_{X\op} * H \\ \hline
\hom_{X\op} * \dd f H \\ \hline
\coend_X \dd f H \\ \hline
\coend_Y H  \\ \hline
\hom_{Y\op} * H 
\end{array}
\]
The fourth one follows from the first one by Proposition \ref{p1}:
\[
\begin{array}{c}
\dm_Z(gf) \\ \hline
\dd\ex_{gf}\hom_X \\ \hline
\dd\ex_g\dd\ex_f\hom_X \\ \hline
\dd\ex_g\hom_Y \\ \hline
\dm_Z g
\end{array}
\]
and implies the fifth one (for $g=\id_Y$) which (again by Proposition \ref{p1}) is equivalent to the first one.
As for the last condition, since a (co)end is a (co)limit weighted by $\hom$
the technique of Section \ref{ip} applies:
\[
\begin{array}{c}
\coend\, g (f\op\tm f) \\ \hline
\hom_{X\op} * g (f\op\tm f)  \\ \hline
\dd\ex_{f\op}\hom_{X\op} * g \\ \hline
\hom_{Y\op} * g \\ \hline
\coend\, g
\end{array}
\qv
\begin{array}{c}
\dd\ex_{f\op}\hom_{X\op} * g \\ \hline
\hom_{X\op} * g (f\op\tm f)  \\ \hline
\coend\, g (f\op\tm f) \\ \hline
\coend\, g  \\ \hline
\hom_{Y\op} * g 
\end{array}
\]
\epf

\begin{remark}
From a two-valued point of view, the second, third and fourth conditions 
of Proposition \ref{ad} become respectively:
\v $H(-,-)$ is reflexive on $Y$ iff $H(f-,f-)$ is so on $X$", 
\v $H(-,-)$ has `fixed points' on $Y$ iff $H(f-,f-)$ has them on $X$" and
\v $g$ and $gf$ have the same (symmetric) image".
\end{remark}

\begin{corollary}
If $f$ and $gf$ are absolutely dense, so it is $g$. 
\epf
\end{corollary}

\begin{prop}
Each of the following is equivalent to the full faithfulness of the functor $f:X\to Y$:
\begin{enumerate}
\item
$\dd f \hom_Y \iso \hom_X$;
\item
$\eend_X H \iso \eend_Y\dd\fa_f H$, naturally in $H$; 
\item
$\coend_X H \iso \coend_Y\dd\ex_f H$, naturally in $H$.
\end{enumerate}
\end{prop}
\pf
The equivalence between the first two conditions is given by
\[
\begin{array}{c}
\eend_Y \dd\fa_f H \\ \hline
\PYY(\hom_Y, \dd\fa_f H) \\ \hline
\PXX(\dd f\hom_Y, H) \\ \hline
\PXX(\hom_X, H) \\ \hline
\eend_X H
\end{array}
\qv
\begin{array}{c}
\PXX(\dd f\hom_Y, H) \\ \hline
\PYY(\hom_Y, \dd\fa_f H) \\ \hline
\eend_Y \dd\fa_f H \\ \hline
\eend_X H  \\ \hline
\PXX(\hom_X, H) 
\end{array}
\]
and similarly one gets the equivalence between the first and the third ones.
\epf


\begin{refs}

\bibitem[El Bashir and Velebil, 2002]{dense} R. El Bashir J. Velebil (2002), Simultaneously Reflective and Coreflective Subcategories of Presheaves, 
{\em Theory and Appl. Cat.} {\bf 10}, 410-423.

\bibitem[Lawvere, 1970]{law70} F.W. Lawvere (1970), {\em Equality in Hyperdoctrines and the Comprehension Scheme as an Adjoint Functor},
Proceedings of the AMS Symposium on Pure Mathematics, XVII, 1-14. 

\bibitem[MacLane, 1965]{mac65} S. MacLane (1965), Categorical Algebra, 
{\em Bull. Am. Math. Soc.} {\bf 71}, 40-106. 

\bibitem[Pisani, 2007]{pis07} C. Pisani (2007), Components, Complements and the Reflection Formula, 
{\em Theory and Appl. Cat.} {\bf 19}, 19-40. 

\bibitem[Pisani, 2008]{pis08} C. Pisani (2008), Balanced Category Theory, 
{\em Theory and Appl. Cat.} {\bf 20}, 85-115. 

\bibitem[Pisani, 2010]{pis10} C. Pisani (2010), A Logic for Categories, 
{\em Theory and Appl. Cat.} {\bf 24}, 394-417.

\bibitem[Wood, 1982]{wood} R. J. Wood (1982), Abstract Proarrows I, 
{\em Cahiers Top. Geo. Diff. Categoriques } {\bf 23(3)}, 279-290.

\end{refs}

\end{document}